\documentclass[11pt,a4paper]{article}
\usepackage[utf8]{inputenc}
\usepackage[T1]{fontenc}
\usepackage{lmodern}
\usepackage[margin=2.5cm]{geometry}
\usepackage{amsmath,amssymb,amsthm}
\usepackage{graphicx,booktabs,xcolor}
\usepackage[numbers,sort&compress]{natbib}
\usepackage{tikz}
\usetikzlibrary{arrows.meta,positioning,shapes,calc,fit}
\usepackage{hyperref}
\hypersetup{colorlinks=true,linkcolor=blue!50!black,citecolor=blue!50!black,urlcolor=blue!50!black}
\newtheorem{theorem}{Theorem}[section]
\newtheorem{proposition}[theorem]{Proposition}
\newtheorem{lemma}[theorem]{Lemma}
\newtheorem{corollary}[theorem]{Corollary}
\theoremstyle{definition}
\newtheorem{assumption}[theorem]{Assumption}
\theoremstyle{remark}
\newtheorem{remark}[theorem]{Remark}

\newcommand{\X}{\mathcal{X}}
\newcommand{\U}{\mathcal{U}}
\newcommand{\R}{\mathbb{R}}
\newcommand{\E}{\mathbb{E}}
\newcommand{\Pb}{\mathbb{P}}

\title{Optimal One Health control of antimicrobial resistance:\\[4pt]
\large a household mean-field model with community feedback and parameters from Burkina Faso}

\author{Mifiamba Soma$^{1,*}$, Pingawind\'e Sawadogo$^{1}$, Yisso Fid\`ele Bacye$^{1}$, Evariste Bako$^{1}$,\\
Mari\'etou Sissao$^{1}$, Maurice Sarigda$^{1}$, Adama Kiemde$^{1}$\\[8pt]
\small $^{1}$Université Thomas SANKARA, Centre Universitaire de Tenkodogo,\\ \small 12 BP 417 Ouagadougou 12, Burkina Faso\\
\small $^{*}$Corresponding author: mifiambasoma@yahoo.fr}

\date{}

\begin{document}
\maketitle

\begin{abstract}
We study how to combine three One Health interventions against antimicrobial resistance: prescription-only dispensing of antibiotics to people, access to veterinary care for animals, and treatment of household sludge and manure. Each household is a Markov chain on the carriage states of its people, animals and soil, and households are linked through the community, whose prevalence feeds the acquisition rates. We prove that, as the number of households grows, the distribution of household states converges to the solution of a nonlinear mean-field equation. This system is monotone for the stochastic order, so that no intervention increases carriage anywhere, and for every constant intervention it has exactly one equilibrium, which attracts all solutions and is given by a scalar equation. For the optimal control problem we prove existence, necessary conditions, positivity of all controls under a quadratic cost, uniqueness for short horizons, and a turnpike property under an explicit dissipativity condition, which we verify for our parameters by a computer-assisted numerical certification. With parameters from Tenkodogo, Burkina Faso, the optimal combined strategy averts about $17\%$ of human carriage over three years in the baseline scenario, and $11\%$ to $28\%$ over the plausible scenarios for the parameters that the survey cannot identify; these figures are projections conditional on the model. The strategy itself is robust: the value of perfect information on the unmeasured parameters is about $2\%$ of the achievable benefit. A compartmental model with the same transmission rates but without households underestimates the benefit, and targeting the interventions by zone of the town increases it.
\end{abstract}



%======================================================================
\section{Introduction}\label{sec:intro}
%======================================================================

Antimicrobial resistance is a major threat to global health \citep{amr2022}, and the global action plan of the World Health Organization calls for coordinated action on human health, animal health and the environment \citep{who2015}. In low-income settings, people, animals and the soil of the same household share resistant bacteria, and they are also exposed to the wider community through markets, water, informal sellers of antibiotics and the spreading of manure \citep{booton2021,chilanga2025}. Public health services can act on several of these routes, but with limited resources, and they need to know which interventions to combine, how strongly, and for how long.

Mathematical models of antimicrobial resistance in a One Health setting have mostly been compartmental models of whole populations \citep{booton2021}. They ignore the household, while the household is where transmission between people and animals takes place and where field studies sample \citep{haverkate2017}. Models of resistance in human populations \citep{bonhoeffer1997,lipsitch2000} and of the effect of reducing antibiotic use in animals \citep{vanbunnik2017,tang2017} rarely represent the household. A household can be described by a continuous-time Markov chain on the carriage states of its people, its animals and the soil of its yard, and the rates of such a chain can be estimated from a household survey; we do this for a survey of $59$ households in Tenkodogo, Burkina Faso (Appendix~\ref{app:estimation}). A single household with fixed acquisition rates from outside is not enough, however. When an intervention lowers carriage in many households, it also lowers the pressure of the community on each household, and a model with fixed external rates misses this feedback.

In this paper we link the households through the community and study the optimal control of the resulting system. We follow the approach of household models in epidemiology \citep{ball1997,house2008}: the state of a household is described exactly, and households interact through community averages. Mathematically, we develop and analyse a class of monotone mean-field optimal control systems generated by interacting household Markov chains (Remark~\ref{rem:general}): we establish a quantitative mean-field limit, uniqueness of the nonlinear equilibrium through a scalar reduction, structural properties of the optimal controls, and a turnpike result based on dissipativity. In more detail, our contributions are the following.
\begin{enumerate}
\item We prove that, as the number $N$ of households grows, the empirical distribution of household states converges to the solution of a nonlinear eight-dimensional mean-field equation, with a bound of order $N^{-1/2}$ that is uniform over the interventions (Theorem~\ref{thm:mf}). An optimal control of the mean-field equation is therefore asymptotically optimal, as $N\to\infty$, for the finite community (Corollary~\ref{cor:mfcost}).
\item We prove that the mean-field system is monotone for the stochastic order of household states (Theorem~\ref{thm:comparison}). In practical terms, no intervention, on any node, can increase carriage on any node at any time, and under constant interventions the prevalences decrease to an equilibrium that depends monotonically on the interventions; all equilibria lie between a smallest and a largest one, which can be computed (Corollary~\ref{cor:norebound} and Propositions~\ref{prop:equilibrium} and~\ref{prop:minmax}).
\item For the optimal control problem with three interventions we prove existence (Theorem~\ref{thm:existence}), derive the optimality system (Theorem~\ref{thm:pmp}), show that the adjoint is increasing in the household state, so that, with a quadratic cost of effort, every intervention is used on $(0,T)$ (Theorem~\ref{thm:positive}), and prove uniqueness for short horizons, together with the convergence of the forward--backward sweep (Theorem~\ref{thm:unique}).
\item We show that the equilibria of the eight-dimensional system are the roots of a scalar equation (Theorem~\ref{thm:scalar}), and, using a sub-homogeneity property of the household chain proved by a graphical construction, that this equation has exactly one root for every constant intervention, so that the equilibrium is unique and globally attracting whatever the rates (Theorem~\ref{thm:unique-eq}); and that the optimal value is monotone in the rates of the model, taken as independent parameters (Corollary~\ref{cor:value}).
\item We give a turnpike theorem: under a strict dissipativity condition, the optimal control on any horizon stays close to the long-run optimum except during a time that does not depend on the horizon, and the constant long-run optimum loses at most a fixed amount (Theorem~\ref{thm:turnpike}). We construct a storage function from the multiplier of the static problem and a Lyapunov function of the linearised dynamics, prove strict dissipativity under an explicit condition (Theorem~\ref{thm:dissip}), and verify it for our parameters by a computer-assisted numerical certification with sums of squares. We also prove an exponential turnpike for the linearised problem, with a rate given by the spectrum of a Hamiltonian matrix (Theorem~\ref{thm:lqturnpike}), which agrees with the rate observed for the nonlinear problem.
\item With parameters from Tenkodogo we compute the optimal strategy, compare all combinations of interventions and check the mean-field limit by exact stochastic simulation. We check the conditions of the turnpike theorems, verify strict dissipativity numerically by sums of squares, compare our model with a compartmental model without households, and compare targeted and uniform strategies in a two-zone model built from the survey. We compute a strategy that is robust to the uncertain parameters and the expected value of perfect information, and we give a cost--effectiveness frontier, a sensitivity study, a robustness study, one- and two-dimensional studies and three-dimensional representations of the optimal controls.
\end{enumerate}

Section~\ref{sec:model} defines the model. Section~\ref{sec:mf} proves the mean-field limit and Section~\ref{sec:monotone} the monotonicity. Section~\ref{sec:control} studies the optimal control problem. Sections~\ref{sec:params} and~\ref{sec:numerics} give the parameters and the numerical results, Section~\ref{sec:guide} a short guide for practitioners, and Section~\ref{sec:discussion} discusses the limits of the study.

%======================================================================
\section{The model}\label{sec:model}
%======================================================================

\subsection{One household}

A household has three nodes: its people ($H$), its animals ($A$) and the soil of its yard ($E$). Each node either carries the resistance trait or not, and the state of the household is $x=(x_H,x_A,x_E)\in\X=\{0,1\}^3$. We order $\X$ componentwise and write $e_k$ for the $k$-th unit vector. Node $k$ acquires the trait at rate
\begin{equation}\label{eq:lambda}
\lambda_k(x;m,u)=(1-u_k)\,a_k(m)+\sum_{j\neq k}\theta_j(u)\,\beta_{jk}\,x_j\qquad(x_k=0),
\end{equation}
and loses it at rate $\gamma_k$ ($x_k=1$). Here $\beta_{jk}\ge 0$ is the rate of transmission from node $j$ to node $k$ inside the household, $a_k(m)$ is the rate of acquisition from outside the household, which depends on the prevalences $m=(m_H,m_A,m_E)$ in the community, and $u=(u_H,u_A,u_E)$ are the interventions; $\theta_H=\theta_A=1$ and $\theta_E(u)=1-u_E$. Other changes of state (two nodes at once) do not occur.

The acquisition from outside has a part that does not depend on the community, for example the use of antibiotics, and a part that is proportional to a weighted community prevalence:
\begin{equation}\label{eq:a}
a_k(m)=\alpha_k\Big[(1-\varphi)+\varphi\,\frac{\omega\cdot m}{\omega\cdot\bar m}\Big]=\alpha_k^{0}+\eta_k\,\omega\cdot m,
\end{equation}
with $\alpha_k^0=(1-\varphi)\alpha_k$ and $\eta_k=\varphi\alpha_k/(\omega\cdot\bar m)$. The parameter $\varphi\in[0,1)$ is the share of outside acquisition due to the community at baseline, $\omega\in[0,\infty)^3$ weights the three community prevalences, and $\bar m$ is the baseline prevalence defined below. With $\varphi=0$ the households are independent and each follows a Markov chain with fixed outside rates.

\subsection{Interventions}

The three interventions act as follows.
\begin{itemize}
\item $u_H$, prescription-only dispensing of antibiotics, in pharmacies and by informal sellers, reduces by the factor $1-u_H$ the acquisition of people from outside the household.
\item $u_A$, access to veterinary care and prudent use of antibiotics in animals, reduces by the factor $1-u_A$ the acquisition of animals from outside.
\item $u_E$, treatment or composting of household sludge and manure before it is spread, reduces by the factor $1-u_E$ the acquisition of the soil from outside and the transmission from the soil to the people and the animals.
\end{itemize}
The rates $\beta_{HE}$ and $\beta_{AE}$ from people and animals to the soil of their own yard stand for direct faecal contamination, for example by animals kept in the yard; treatment of collected sludge and manure does not prevent it, and we leave these rates unchanged. All results of Sections~\ref{sec:mf}--\ref{sec:control} remain true if $u_E$ also reduces them, since they only use the fact that the rates are nonincreasing in $u$.
An intervention is a measurable function $u:[0,T]\to[0,\bar u]^3$, with $\bar u<1$; we write $\U$ for the set of interventions. The value $u_k(t)$ is the relative reduction achieved at time $t$ and $\bar u$ expresses that no intervention removes a route completely. Measured effects of such measures are scarce and variable. In a meta-analysis of interventions restricting the use of antibiotics in food animals, the absolute reduction of resistant bacteria in the animals was commonly $10$ to $15\%$, and the prevalence of resistance in people was $24\%$ lower in the intervention groups \citep{tang2017}. In a cluster randomised trial in rural Burkina Faso, a community bundle on antibiotic use and hygiene reduced the rate of household acquisition of ESBL-producing \emph{E.~coli} by an estimated $12\%$, with a credible interval from a reduction of $35\%$ to an increase of $23\%$ \citep{aizouk2025}. The optimal values of $u_k$ found in Section~\ref{sec:numerics}, mostly between $0.1$ and $0.3$, are of this order.

\subsection{Many households and the mean-field equation}

Consider $N$ households with states $X^1(t),\dots,X^N(t)$. Household $i$ changes state with the rates~\eqref{eq:lambda}, in which $m$ is the empirical prevalence
\[
m^N(t)=\frac1N\sum_{i=1}^N X^i(t),
\]
and $u$ is the same for all households. We write $P^N(t)$ for the empirical distribution of household states, a probability vector on $\X$, so that $m^N=P^NS$, where $S$ is the $8\times 3$ matrix whose rows are the states.

For $m\in[0,1]^3$ and $u\in[0,\bar u]^3$ let $Q(m,u)$ be the $8\times 8$ generator of the household chain with rates~\eqref{eq:lambda}. Let $\Delta$ be the set of probability (row) vectors on $\X$. The mean-field equation is
\begin{equation}\label{eq:mf}
\dot p(t)=p(t)\,Q\big(p(t)S,\,u(t)\big),\qquad p(0)=p_0\in\Delta.
\end{equation}
It is a system of eight nonlinear ordinary differential equations; each component $p(t,x)$ is the proportion of households in state $x$. The prevalences are $m(t)=p(t)S$. The household structure is kept exactly: the model is not a compartmental model of three populations, since the joint state of the three nodes of a household matters for the transmission inside it.

\emph{Baseline.} Let $\bar p$ be the stationary law of the household chain with outside rates $\alpha$ and no intervention, that is $\bar pQ_0=0$ with $Q_0$ the generator with rates $\alpha_k+\sum_j\beta_{jk}x_j$, and let $\bar m=\bar pS$. By~\eqref{eq:a}, $a(\bar m)=\alpha$, so $\bar p$ is an equilibrium of~\eqref{eq:mf} with $u=0$: the parameters $\alpha$ estimated from a survey at one time (Appendix~\ref{app:estimation}) are the outside rates at the present equilibrium of the community.

\begin{assumption}\label{ass:main}
$\alpha_k>0$, $\gamma_k>0$ and $\beta_{jk}\ge 0$ for all nodes; $\varphi\in[0,1)$; $\omega\in[0,\infty)^3$ with $\omega\cdot\bar m>0$; and $\bar u\in(0,1)$.
\end{assumption}

\begin{remark}[General framework]\label{rem:general}
Nothing in Sections~\ref{sec:mf}--\ref{sec:control} uses the fact that a household has three nodes. The results hold for households with $n$ nodes, state space $\{0,1\}^n$, acquisition rates $(1-u_k)a_k(m)+\sum_{j\neq k}\theta_j(u)\beta_{jk}x_j$ with $a_k(m)=\alpha_k^0+\eta_k\,\omega\cdot m$, $\alpha^0_k>0$, $\eta_k,\omega_j,\beta_{jk}\ge0$, interventions that only reduce rates, and loss rates $\gamma_k>0$: that is, for mean-field limits of interacting household Markov chains that are attractive, additive and coupled through a scalar community pressure. The constants depend on the number $2^n$ of household states (for example $\sqrt8$ becomes $2^{n/2}$ in Theorem~\ref{thm:mf}, and the linearised problem of Section~\ref{sec:control} has dimension $2^n-1$). The three-node model of the One Health application is one instance of this class; several people and animals per household, or water and food as further nodes, are others.
\end{remark}

We use the norm $\|p\|_1=\sum_x|p(x)|$ on $\R^8$ and write
\[
\bar q=\sum_{k}\Big(\alpha_k^0+\eta_k\textstyle\sum_j\omega_j+\sum_{j\neq k}\beta_{jk}+\gamma_k\Big),\qquad \hat\eta=\sum_k\eta_k,\qquad\hat\omega=\sum_j\omega_j .
\]
The total rate of change of a household is at most $\bar q$, for every state, prevalence and intervention.

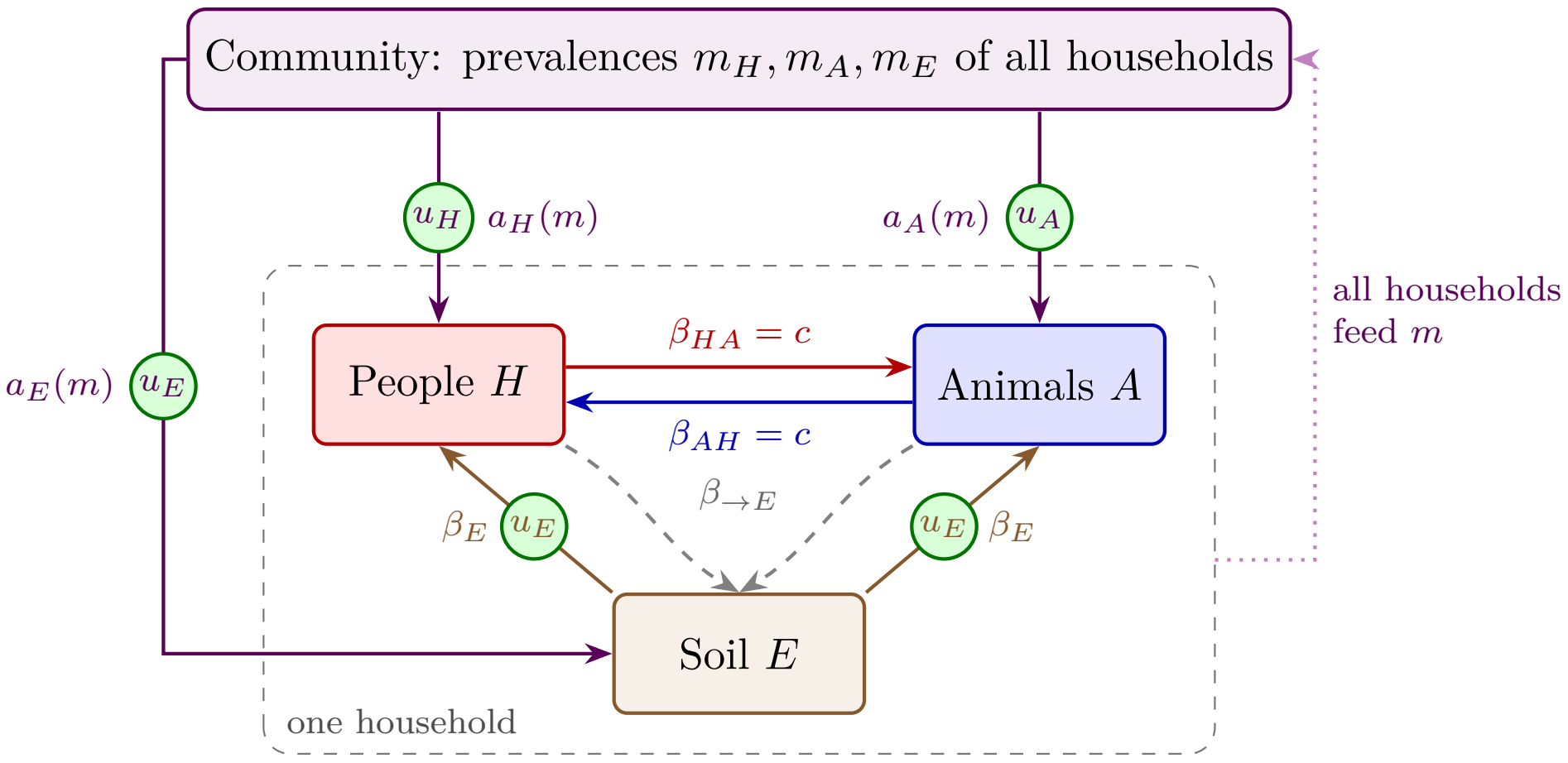
\begin{figure}[t]
\centering
\begin{tikzpicture}[
  hh/.style={draw=#1!70!black, fill=#1!12, thick, rounded corners=3pt, minimum width=2.0cm, minimum height=0.95cm, align=center, font=\small},
  arr/.style={-{Stealth[length=2.2mm]}, thick},
  lab/.style={font=\scriptsize},
  valve/.style={circle, draw=green!45!black, fill=green!15, thick, inner sep=1pt, font=\scriptsize}]
\draw[black!55, dashed, rounded corners=6pt] (-3.8,-2.15) rectangle (3.8,1.75);
\node[font=\scriptsize, black!70, anchor=south west] at (-3.75,-2.12) {one household};
\node[hh=red] (H) at (-2.4,0.8) {People $H$};
\node[hh=blue] (A) at (2.4,0.8) {Animals $A$};
\node[hh=brown] (Ev) at (0,-1.35) {Soil $E$};
\node[draw=violet!70!black, fill=violet!8, thick, rounded corners=4pt, minimum width=8.4cm, minimum height=0.8cm, font=\small] (C) at (0,3.4) {Community: prevalences $m_H, m_A, m_E$ of all households};
% within household
\draw[arr, red!70!black] ([yshift=4pt]H.east) -- node[above, lab] {$\beta_{HA}=c$} ([yshift=4pt]A.west);
\draw[arr, blue!70!black] ([yshift=-4pt]A.west) -- node[below, lab] {$\beta_{AH}=c$} ([yshift=-4pt]H.east);
\draw[arr, brown!70!black] (Ev.north west) -- node[pos=0.45, valve] {$u_E$} node[pos=0.45, left=6pt, lab] {$\beta_E$} (H.south);
\draw[arr, brown!70!black] (Ev.north east) -- node[pos=0.45, valve] {$u_E$} node[pos=0.45, right=6pt, lab] {$\beta_E$} (A.south);
\draw[arr, gray, dashed] (H.south east) to[out=-30, in=150] (Ev.north);
\draw[arr, gray, dashed] (A.south west) to[out=-150, in=30] (Ev.north);
\node[lab, gray!80!black] at (0,-0.08) {$\beta_{\to E}$};
% community
\draw[arr, violet!70!black] (-2.4,2.98) -- node[pos=0.5, valve] {$u_H$} node[pos=0.5, right=7pt, lab] {$a_H(m)$} (H.north);
\draw[arr, violet!70!black] (2.4,2.98) -- node[pos=0.5, valve] {$u_A$} node[pos=0.5, left=7pt, lab] {$a_A(m)$} (A.north);
\draw[arr, violet!70!black] (C.west) -- (-4.6,3.4) -- node[pos=0.55, valve] {$u_E$} node[pos=0.55, left=7pt, lab] {$a_E(m)$} (-4.6,-1.35) -- (Ev.west);
% feedback
\draw[arr, violet!50, dotted] (3.8,-0.6) -- (4.6,-0.6) -- node[right, lab, align=left, violet!70!black] {all households\\feed $m$} (4.6,3.4) -- (C.east);
\end{tikzpicture}
\caption{The model. Each household is a Markov chain on the carriage states of its people, animals and soil, with transmission inside the household (arrows between nodes) and acquisition from the community (violet arrows), whose rate depends on the prevalences in all households. Each node also loses the trait at rate $\gamma_k$ (not drawn). The green circles mark the routes reduced by the three interventions: $u_H$, prescription-only dispensing; $u_A$, access to veterinary care; $u_E$, treatment of sludge and manure. A route marked $u_k$ has its rate multiplied by $1-u_k$.}
\label{fig:model}
\end{figure}

%======================================================================
\section{Well-posedness and mean-field limit}\label{sec:mf}
%======================================================================

\begin{proposition}\label{prop:wellposed}
Under Assumption~\ref{ass:main}, for every $u\in\U$ and $p_0\in\Delta$, equation~\eqref{eq:mf} has a unique absolutely continuous solution on $[0,T]$, and $p(t)\in\Delta$ for all $t$. If $t>0$, then $p(t,x)>0$ for every $x\in\X$. For $p,p'\in\Delta$ and every $u$,
\begin{equation}\label{eq:lip}
\big\|pQ(pS,u)-p'Q(p'S,u)\big\|_1\le L\,\|p-p'\|_1,\qquad L=2(\bar q+\hat\eta\hat\omega),
\end{equation}
and the map $u\mapsto p$ is Lipschitz from $L^1(0,T)^3$ to $C([0,T])^8$.
\end{proposition}

\begin{proof}
To work on an open set, we extend $Q(m,u)$ to all $m\in\R^3$ by replacing $m$ with its projection on $[0,1]^3$; the extension is Lipschitz, and it coincides with $Q$ on $\Delta S\subset[0,1]^3$, so it does not change the solutions in $\Delta$. The right-hand side $F(p,t)=pQ(pS,u(t))$ is then locally Lipschitz in $p$, with coefficients that are bounded measurable functions of $t$, so local existence and uniqueness of an absolutely continuous solution follow from Carathéodory's theorem. The rows of $Q$ sum to zero, so $\sum_xp(t,x)$ is constant. Each component satisfies $\dot p(t,x)\ge -q(x,t)p(t,x)$, where $q(x,t)\le\bar q$ is the total exit rate of $x$ at time $t$, as long as $p\ge 0$; hence $p(t,x)\ge p(0,x)e^{-\bar qt}\ge 0$, and the solution stays in the compact set $\Delta$ and exists on $[0,T]$. For positivity, the household chain with generator $Q(p(t)S,u(t))$ has acquisition rates at least $(1-\bar u)\alpha_k^0>0$ and loss rates $\gamma_k>0$, so any state can be reached from any other state in at most three steps with rates bounded below; the law $p(t)$ of this time-inhomogeneous chain therefore charges every state for $t>0$.

For~\eqref{eq:lip}, write $m=pS$, $m'=p'S$ and
\[
pQ(m,u)-p'Q(m',u)=(p-p')Q(m,u)+p'\big(Q(m,u)-Q(m',u)\big).
\]
Since every row of $Q(m,u)$ has absolute sum at most $2\bar q$, the first term has norm at most $2\bar q\|p-p'\|_1$. The rates~\eqref{eq:lambda} depend on $m$ only through $a_k(m)$, and $|a_k(m)-a_k(m')|\le\eta_k\hat\omega\max_j|m_j-m'_j|\le\eta_k\hat\omega\|p-p'\|_1$, so every row of $Q(m,u)-Q(m',u)$ has absolute sum at most $2\hat\eta\hat\omega\|p-p'\|_1$, and $p'\in\Delta$ gives the second bound. Finally $Q$ is affine in $u$ with bounded coefficients, so $\|F(p,u)-F(p,v)\|_1\le K|u-v|$ for a constant $K$, and Gronwall's lemma gives $\|p_u(t)-p_v(t)\|_1\le Ke^{LT}\int_0^t|u-v|\,ds$.
\end{proof}

\begin{theorem}[Mean-field limit]\label{thm:mf}
Let Assumption~\ref{ass:main} hold, let $u\in\U$, and let $p$ be the solution of~\eqref{eq:mf}. Then
\begin{equation}\label{eq:mfbound}
\E\Big[\sup_{0\le t\le T}\|P^N(t)-p(t)\|_1\Big]\le\Big(\E\|P^N(0)-p_0\|_1+8\sqrt{\bar qT/N}\Big)e^{LT},
\end{equation}
and, for every $t\in[0,T]$,
\begin{equation}\label{eq:mfbound2}
\E\|P^N(t)-p(t)\|_1\le\Big(\E\|P^N(0)-p_0\|_1+4\sqrt{2\bar qt/N}\Big)e^{2\hat\eta\hat\omega t}.
\end{equation}
If the initial states of the households are independent with law $p_0$, then $\E\|P^N(0)-p_0\|_1\le\sqrt{8/N}$. The bounds hold for every measurable $u\in\U$, with constants that depend on $T$ and on the rates through $\bar q$, $L$ and $\hat\eta\hat\omega$, but not on $u$ or on $\bar u$, since $\bar q$ bounds the rates without intervention. Here $\|\cdot\|_1$ is the $\ell^1$ norm on $\R^8$ and the expectation is over the random evolution of the $N$ households.
\end{theorem}

\begin{proof}
The process $P^N$ is a time-inhomogeneous Markov chain on $\Delta$. When one household moves from $x$ to $y$, $P^N$ changes by $(\delta_y-\delta_x)/N$, and this happens at rate $NP^N(x)\,q(x,y;m^N,u(t))$, where $q(x,y;m,u)$ are the entries of $Q(m,u)$. The drift of $P^N$ is therefore
\[
\sum_{x\neq y}NP^N(x)\,q(x,y;m^N,u)\,\frac{\delta_y-\delta_x}{N}=P^NQ(P^NS,u),
\]
and
\[
M^N(t)=P^N(t)-P^N(0)-\int_0^tP^N(s)Q\big(P^N(s)S,u(s)\big)\,ds
\]
is a martingale. Each jump has Euclidean norm $\sqrt2/N$, and the total jump rate is at most $N\bar q$, so the trace of the predictable quadratic variation of $M^N$ is at most $2\bar qt/N$. By Doob's inequality, $\E\sup_{t\le T}\|M^N(t)\|_2^2\le 8\bar qT/N$, and since $\|v\|_1\le\sqrt8\|v\|_2$ on $\R^8$, $\E\sup_{t\le T}\|M^N(t)\|_1\le 8\sqrt{\bar qT/N}$. Subtracting the integral form of~\eqref{eq:mf} and using~\eqref{eq:lip},
\[
\|P^N(t)-p(t)\|_1\le\|P^N(0)-p_0\|_1+\sup_{s\le T}\|M^N(s)\|_1+L\int_0^t\|P^N(s)-p(s)\|_1\,ds,
\]
and Gronwall's lemma gives~\eqref{eq:mfbound}.

For~\eqref{eq:mfbound2}, let $\Pi(s,t)$ be the propagator of the linear equation $\dot\pi=\pi Q(m(t),u(t))$, with $m$ the mean-field prevalence. Each $\Pi(s,t)$ is a stochastic matrix, so $\|v\Pi(s,t)\|_1\le\|v\|_1$. The error $e=P^N-p$ satisfies $de=\big[eQ(m,u)+P^N\big(Q(m^N,u)-Q(m,u)\big)\big]dt+dM^N$, and by the variation of constants formula, which holds path by path because $M^N$ has finite variation,
\[
e(t)=e(0)\Pi(0,t)+\int_0^tP^N(s)\big(Q(m^N(s),u(s))-Q(m(s),u(s))\big)\Pi(s,t)\,ds+\int_0^tdM^N(s)\,\Pi(s,t).
\]
The norm of the first term is at most $\|e(0)\|_1$, and that of the integrand of the second at most $2\hat\eta\hat\omega\|e(s)\|_1$, as in the proof of~\eqref{eq:lip}. For fixed $t$, $r\mapsto\int_0^rdM^N(s)\Pi(s,t)$ is a martingale whose jumps $(\delta_y-\delta_x)\Pi(s,t)/N$ have Euclidean norm at most their $\ell^1$ norm, that is at most $2/N$; with a total jump rate at most $N\bar q$, its second moment at $r=t$ is at most $4\bar qt/N$, and its expected $\ell^1$ norm at most $\sqrt8\sqrt{4\bar qt/N}=4\sqrt{2\bar qt/N}$. Gronwall's lemma applied to $\E\|e(t)\|_1$ gives~\eqref{eq:mfbound2}. For independent initial states, $\E\|P^N(0)-p_0\|_2^2=\sum_xp_0(x)(1-p_0(x))/N\le 1/N$, so $\E\|P^N(0)-p_0\|_1\le\sqrt{8/N}$.
\end{proof}

The proof is the classical argument of \citet{kurtz1970} \citep[see also][Chapter~11]{ethier1986}; we give it because we need a bound that is explicit and uniform in the intervention. The result is a law of large numbers for a system of households in which each household interacts with the others only through the community prevalences; by exchangeability, any fixed number of households become asymptotically independent, each with law $p(t)$.

The bound~\eqref{eq:mfbound2} grows with the strength $\hat\eta\hat\omega$ of the community feedback only; without community ($\varphi=0$, so $\eta=0$) it has no exponential factor, as expected for independent households. The constants are worst-case bounds: with the parameters of Section~\ref{sec:params}, $\bar q\approx7.3$ and $2\hat\eta\hat\omega\approx2.4$, so for a horizon of ten time units the factor $e^{2\hat\eta\hat\omega T}$ is too large for the bound to be informative at realistic sizes, and the rate $N^{-1/2}$ is the useful part of the theorem. Section~\ref{sec:numerics} shows by simulation that the actual error is small from a few thousand households.

For a finite community the natural cost of an intervention is
\[
J^N(u)=\E\int_0^T\Big[w\cdot m^N(t)+\tfrac12\sum_kB_ku_k(t)^2\Big]dt,
\]
where $w\in[0,\infty)^3$ weights the prevalences and $B_k>0$ is the cost of intervention $k$; $J(u)$ is the same quantity for the solution of~\eqref{eq:mf} (Section~\ref{sec:control}).

\begin{corollary}\label{cor:mfcost}
Under the assumptions of Theorem~\ref{thm:mf} with independent initial states, $\sup_{u\in\U}|J^N(u)-J(u)|\le\delta_N$, with $\delta_N=T(w_H+w_A+w_E)\big(\sqrt8+4\sqrt{2\bar qT}\big)e^{2\hat\eta\hat\omega T}N^{-1/2}$. Hence, if $u^*$ minimises $J$ over $\U$, then $J^N(u^*)\le\inf_{u\in\U}J^N(u)+2\delta_N$: among interventions fixed in advance, an optimal intervention for the mean-field equation is asymptotically optimal for the community of $N$ households as $N\to\infty$.
\end{corollary}

\begin{proof}
$|w\cdot(m^N-m)|\le\sum_kw_k|m^N_k-m_k|\le(\sum_kw_k)\|P^N-p\|_1$, and~\eqref{eq:mfbound2}, integrated over $[0,T]$, gives the first statement. Then $J^N(u^*)\le J(u^*)+\delta_N\le J(u)+\delta_N\le J^N(u)+2\delta_N$ for every $u$.
\end{proof}

%======================================================================
\section{Monotonicity}\label{sec:monotone}
%======================================================================

A set $\Gamma\subset\X$ is increasing if $x\in\Gamma$ and $x\le x'$ imply $x'\in\Gamma$. For $p,p'\in\Delta$ we write $p\preceq p'$ if $p(\Gamma)\le p'(\Gamma)$ for every increasing set $\Gamma$ (stochastic order). Since $\{x:x_k=1\}$ is increasing, $p\preceq p'$ implies $pS\le p'S$ componentwise. By Strassen's theorem \citep{strassen1965,kamae1977}, $p\preceq p'$ if and only if there are random states $X\le X'$ with laws $p$ and $p'$.

The household chain is attractive: by~\eqref{eq:lambda}, the acquisition rate of a node is a nondecreasing function of the states of the other nodes, and the loss rates do not depend on the state. The following lemma is a time-inhomogeneous version of classical comparison results for such chains \citep{kamae1977,massey1987,liggett1985}; we give a short proof.

\begin{lemma}\label{lem:coupling}
Consider two single-flip chains on $\X$ with acquisition rates $\lambda_k(x,t)$ and $\tilde\lambda_k(x,t)$, and loss rates $\gamma_k$ and $\tilde\gamma_k$ not depending on the state, all nonnegative, measurable and bounded in $t$. Assume that $\gamma_k\ge\tilde\gamma_k$ and that for every $t$, every node $k$ and all $x\le\tilde x$ with $x_k=\tilde x_k=0$, $\lambda_k(x,t)\le\tilde\lambda_k(\tilde x,t)$. If the initial laws satisfy $\pi_0\preceq\tilde\pi_0$, then the laws at time $t$ satisfy $\pi(t)\preceq\tilde\pi(t)$ for all $t\in[0,T]$.
\end{lemma}

\begin{proof}
Let $\bar r$ bound all the rates. Take $X(0)\le\tilde X(0)$ with laws $\pi_0$ and $\tilde\pi_0$, and, independently, for each node $k$ two Poisson processes of rate $\bar r$ with independent uniform marks $V$ on $[0,1]$. At an event of the first process at time $s$ with mark $V$, node $k$ of $X$ acquires the trait if $X_k(s^-)=0$ and $V\bar r<\lambda_k(X(s^-),s)$, and similarly for $\tilde X$ with $\tilde\lambda_k$. At an event of the second process, node $k$ of $X$ loses the trait if it carries it and $V\bar r<\gamma_k$, and node $k$ of $\tilde X$ if it carries it and $V\bar r<\tilde\gamma_k$. By thinning, $X$ and $\tilde X$ are Markov chains with the given rates. The order $X\le\tilde X$ is preserved at every event: at a loss event with $X_k=\tilde X_k=1$, $\tilde\gamma_k\le\gamma_k$ implies that $X$ loses the trait whenever $\tilde X$ does; at an acquisition event with $X_k=\tilde X_k=0$, the assumption gives $\lambda_k(X,s)\le\tilde\lambda_k(\tilde X,s)$, so $\tilde X$ acquires whenever $X$ does; and if $X_k=0$, $\tilde X_k=1$ the order cannot be broken. Hence $X(t)\le\tilde X(t)$ for all $t$.
\end{proof}

\begin{theorem}[Comparison]\label{thm:comparison}
Let $u,\tilde u\in\U$ with $u(t)\ge\tilde u(t)$ componentwise for almost every $t$, and let $p_0\preceq\tilde p_0$. Let $p$ and $\tilde p$ be the solutions of~\eqref{eq:mf} with $(u,p_0)$ and $(\tilde u,\tilde p_0)$. Then $p(t)\preceq\tilde p(t)$, and in particular $m(t)\le\tilde m(t)$, for all $t\in[0,T]$.
\end{theorem}

\begin{proof}
For a continuous path $\mu:[0,T]\to[0,1]^3$, let $\Phi_{u}(\mu)$ be the law, at all times, of the linear chain with rates $\lambda_k(x;\mu(t),u(t))$ and initial law $p_0$, that is the solution of $\dot\pi=\pi Q(\mu(t),u(t))$, $\pi(0)=p_0$; define $\tilde\Phi_{\tilde u}$ in the same way with $\tilde p_0$. A function $p$ solves~\eqref{eq:mf} if and only if $p=\Phi_u(pS)$.

\emph{Step 1: comparison of the linear chains.} If $\mu\le\tilde\mu$, the rates~\eqref{eq:lambda} satisfy, for $x\le\tilde x$ with $x_k=\tilde x_k=0$,
\[
\lambda_k(x;\mu,u)=(1-u_k)a_k(\mu)+\sum_j\theta_j(u)\beta_{jk}x_j\le(1-\tilde u_k)a_k(\tilde\mu)+\sum_j\theta_j(\tilde u)\beta_{jk}\tilde x_j=\lambda_k(\tilde x;\tilde\mu,\tilde u),
\]
because $a_k$ is nondecreasing, $\theta_j$ and $1-u_k$ are nonincreasing in $u$, and $x\le\tilde x$. Lemma~\ref{lem:coupling} gives $\Phi_u(\mu)(t)\preceq\tilde\Phi_{\tilde u}(\tilde\mu)(t)$, and in particular $\Phi_u(\mu)S\le\tilde\Phi_{\tilde u}(\tilde\mu)S$.

\emph{Step 2: iteration.} Let $\mu^0=m=pS$ and $\mu^{n+1}=\tilde\Phi_{\tilde u}(\mu^n)S$. By Step~1 with $\mu=\tilde\mu=m$, $\mu^1\ge\Phi_u(m)S=m=\mu^0$, and by induction $\mu^{n+1}=\tilde\Phi_{\tilde u}(\mu^n)S\ge\tilde\Phi_{\tilde u}(\mu^{n-1})S=\mu^n$ (Step~1 with $u=\tilde u$ and the same initial law). As in the proof of~\eqref{eq:lip}, two linear chains driven by $\mu$ and $\mu'$ satisfy $\|\tilde\Phi(\mu)(t)-\tilde\Phi(\mu')(t)\|_1\le K\int_0^t\max_j|\mu_j-\mu'_j|\,ds$ for a constant $K$, so $\sup_{s\le t}|\mu^{n+1}-\mu^n|\le(Kt)^n/n!\,\sup|\mu^1-\mu^0|$, and $\mu^n$ converges uniformly to a fixed point of $\mu\mapsto\tilde\Phi_{\tilde u}(\mu)S$. The fixed point gives a solution of~\eqref{eq:mf} with $(\tilde u,\tilde p_0)$, which is unique by Proposition~\ref{prop:wellposed}, so $\mu^n\to\tilde m$ and $m\le\tilde m$. Step~1 with $\mu=m$ and $\tilde\mu=\tilde m$ then gives $p(t)=\Phi_u(m)(t)\preceq\tilde\Phi_{\tilde u}(\tilde m)(t)=\tilde p(t)$.
\end{proof}

\begin{corollary}[No rebound]\label{cor:norebound}
Start from the baseline equilibrium $\bar p$. For every intervention $u\in\U$, $m_k(t)\le\bar m_k$ for every node $k$ and every $t$. If $u\ge\tilde u$, then $m(t)\le\tilde m(t)$ for all $t$. Thus, in this model, an intervention on any node never increases carriage on any node at any time, and more intervention never gives more carriage.
\end{corollary}

\begin{proof}
Apply Theorem~\ref{thm:comparison} with $\tilde u=0$ and $p_0=\tilde p_0=\bar p$, whose solution is constant, and then with $\tilde u$ in place of $0$.
\end{proof}

The proof of Theorem~\ref{thm:comparison} uses only the order of the rates, and Lemma~\ref{lem:coupling} allows different loss rates. The same argument therefore gives the dependence on the parameters. Write the community term as $a_k(m)=\alpha_k^0+\eta_k\,\omega\cdot m$.

\begin{proposition}[Dependence on the parameters]\label{prop:params}
Let two sets of parameters satisfy $\alpha^0\le\tilde\alpha^0$, $\eta\le\tilde\eta$, $\omega\le\tilde\omega$, $\beta\le\tilde\beta$ and $\gamma\ge\tilde\gamma$ componentwise, and let $u\ge\tilde u$ and $p_0\preceq\tilde p_0$. Then the corresponding solutions of~\eqref{eq:mf} satisfy $p(t)\preceq\tilde p(t)$ for all $t$.
\end{proposition}

\begin{proof}
In Step~1 of the proof of Theorem~\ref{thm:comparison}, the inequality between the acquisition rates still holds, because every term of~\eqref{eq:lambda} is nondecreasing in $\alpha^0,\eta,\omega,\beta$ and in $\mu$; Lemma~\ref{lem:coupling} applies with $\gamma\ge\tilde\gamma$. Step~2 is unchanged.
\end{proof}

In this proposition $\alpha^0$, $\eta$ and $\omega$ are regarded as independent parameters. In the calibrated parametrisation $(\alpha,\varphi,\omega)$ of~\eqref{eq:a}, they are linked through $\bar m$, which itself depends on $\alpha$, $\beta$ and $\gamma$; the proposition then gives no monotonicity in $\varphi$, and none in $\beta$ or $\gamma$ at fixed $\alpha$.

This conclusion does not hold in every model. In a model in which the community acted on the households through competition between strains, or in which an intervention displaced exposure from one route to another, an intervention could increase carriage somewhere. For one household with fixed outside rates, the same coupling argument shows that the stationary prevalences are nondecreasing in every acquisition and transmission rate; Corollary~\ref{cor:norebound} extends this to the community level, where the outside rates themselves respond to the interventions.

\begin{proposition}[Constant interventions]\label{prop:equilibrium}
Let $u\in[0,\bar u]^3$ be constant and let $p$ solve~\eqref{eq:mf} with $p(0)=\bar p$. Then $t\mapsto p(t)$ is nonincreasing for $\preceq$, and $p(t)$ converges, as $t\to\infty$, to an equilibrium $p^\infty(u)$ of~\eqref{eq:mf}. If $u\le v$ componentwise, then $p^\infty(v)\preceq p^\infty(u)$; in particular the long-run prevalences are nonincreasing functions of each intervention.
\end{proposition}

\begin{proof}
Let $h>0$. By Corollary~\ref{cor:norebound}, $p(h)\preceq\bar p=p(0)$. The equation is autonomous, so $t\mapsto p(t+h)$ is the solution started at $p(h)$, and Theorem~\ref{thm:comparison} gives $p(t+h)\preceq p(t)$ for all $t$. Hence $p(t)(\Gamma)$ is nonincreasing for every increasing set $\Gamma$, and converges. The indicator of any single state is a difference of indicators of increasing sets, so $p(t)$ converges to some $p^\infty\in\Delta$. Since the right-hand side of~\eqref{eq:mf} is continuous and $p(t)$ converges, $\dot p(t)$ converges, and its limit must be zero, so $p^\infty$ is an equilibrium. For $u\le v$, Theorem~\ref{thm:comparison} gives $p_v(t)\preceq p_u(t)$ for all $t$, and the order passes to the limit.
\end{proof}

\begin{proposition}[Smallest and largest equilibria]\label{prop:minmax}
Let $u\in[0,\bar u]^3$ be constant. The solution of~\eqref{eq:mf} started from $\delta_{(0,0,0)}$ (no node carries the trait) is nondecreasing for $\preceq$ and converges to an equilibrium $p_{\min}(u)$; the solution started from $\delta_{(1,1,1)}$ is nonincreasing and converges to an equilibrium $p_{\max}(u)$. Every equilibrium $p^e$ satisfies $p_{\min}(u)\preceq p^e\preceq p_{\max}(u)$, and for every initial law the limit points of the solution lie between $p_{\min}(u)$ and $p_{\max}(u)$. In particular, if $p_{\min}(u)=p_{\max}(u)$, the equilibrium is unique and attracts every solution.
\end{proposition}

\begin{proof}
Since $\delta_{(0,0,0)}\preceq p(h)$ for every $h$, the argument of Proposition~\ref{prop:equilibrium}, with the inequalities reversed, shows that the solution from $\delta_{(0,0,0)}$ is nondecreasing and converges to an equilibrium; the solution from $\delta_{(1,1,1)}$ is treated in the same way. If $p^e$ is an equilibrium, or any initial law, then $\delta_{(0,0,0)}\preceq p^e\preceq\delta_{(1,1,1)}$, and Theorem~\ref{thm:comparison} with $\tilde u=u$ orders the three solutions at all times; letting $t\to\infty$ gives the statement.
\end{proof}

The equilibria can be described more precisely. Since the community acts on a household only through the scalar pressure $s=\omega\cdot m$, they reduce to the roots of a scalar equation. For a constant $u$ and $s\ge0$, let $\pi_u(s)$ be the stationary law of the household chain with outside rates $(1-u_k)(\alpha_k^0+\eta_ks)$ and the within-household rates of~\eqref{eq:lambda}; this chain is irreducible, so $\pi_u(s)$ is unique. Let
\[
g_u(s)=\omega\cdot\pi_u(s)S ,
\]
the community pressure that the households would produce if they were exposed to the pressure $s$.

\begin{theorem}[Scalar equation for the equilibria]\label{thm:scalar}
Let $u\in[0,\bar u]^3$ be constant.
\begin{enumerate}
\item[(i)] $p\in\Delta$ is an equilibrium of~\eqref{eq:mf} if and only if $p=\pi_u(s)$ for a root $s$ of $g_u(s)=s$; then $s=\omega\cdot pS$.
\item[(ii)] $g_u$ is a rational function of $s$, nondecreasing on $[0,\infty)$, with $0\le g_u\le\hat\omega$, and it is nonincreasing in each $u_k$. The equation $g_u(s)=s$ has at least one root in $[0,\hat\omega]$. If $s_-$ and $s_+$ are its smallest and largest roots, then $p_{\min}(u)=\pi_u(s_-)$ and $p_{\max}(u)=\pi_u(s_+)$.
\item[(iii)] If $g_u(s)=s$ has only one root, in particular if $g_u'(s)<1$ on $[0,\hat\omega]$, the equilibrium is unique and every solution of~\eqref{eq:mf} with the constant intervention $u$ converges to it.
\end{enumerate}
\end{theorem}

\begin{proof}
(i) If $p$ is an equilibrium, then $pQ(pS,u)=0$, so $p$ is a stationary law of the household chain with outside rates $(1-u_k)a_k(pS)=(1-u_k)(\alpha_k^0+\eta_ks)$, $s=\omega\cdot pS$. By uniqueness, $p=\pi_u(s)$, and $s=\omega\cdot\pi_u(s)S=g_u(s)$. Conversely, if $g_u(s)=s$ and $p=\pi_u(s)$, then $\omega\cdot pS=s$, so the outside rates of $Q(pS,u)$ are those used to define $\pi_u(s)$, and $pQ(pS,u)=0$.

(ii) $\pi_u(s)$ solves a linear system whose coefficients are affine in $s$, so by Cramer's rule it is rational in $s$. If $s\le s'$, the rates with $s'$ dominate those with $s$ in the sense of Lemma~\ref{lem:coupling}; starting both chains from the same law, their laws are ordered at all times, and letting $t\to\infty$ gives $\pi_u(s)\preceq\pi_u(s')$, hence $g_u(s)\le g_u(s')$ since $\omega\ge0$. In the same way $g_u$ is nonincreasing in $u$. Since the prevalences lie in $[0,1]$, $0\le g_u\le\hat\omega$, so $s\mapsto g_u(s)-s$ is nonnegative at $0$ and nonpositive at $\hat\omega$ and has a root. The equilibrium $p_{\min}(u)$ equals $\pi_u(s_*)$ for some root $s_*$ by (i). For every root $s$, Proposition~\ref{prop:minmax} gives $p_{\min}(u)\preceq\pi_u(s)$, so $s_*=\omega\cdot p_{\min}(u)S\le\omega\cdot\pi_u(s)S=s$; hence $s_*=s_-$. The same argument gives $p_{\max}(u)=\pi_u(s_+)$.

(iii) If there is only one root, (ii) gives $p_{\min}(u)=p_{\max}(u)$, and Proposition~\ref{prop:minmax} gives the conclusion. If $g_u'<1$ on $[0,\hat\omega]$, the function $g_u(s)-s$ is strictly decreasing there, and all roots lie in $[0,\hat\omega]$.
\end{proof}

Theorem~\ref{thm:scalar} reduces the study of the equilibria of an eight-dimensional nonlinear system to a scalar equation. We now show that this equation always has exactly one root. The proof rests on a property of the household chain with respect to its outside rates. Let $\pi(a)$ denote the stationary law of the household chain with outside rates $a\in(0,\infty)^3$ and fixed transmission and loss rates, and $m(a)=\pi(a)S$.

\begin{lemma}[Sub-homogeneity]\label{lem:subhom}
For every $a\in(0,\infty)^3$ and $\vartheta\in(0,1)$, $m(\vartheta a)\ge\vartheta\,m(a)$ componentwise.
\end{lemma}

\begin{proof}
The household chain has a graphical representation \citep{liggett1985}: for each node $k$ take independent Poisson processes of \emph{arrivals} at rate $a_k$, of \emph{recovery marks} at rate $\gamma_k$, and, for each $j\neq k$, of \emph{arrows} from $j$ to $k$ at rate $\beta_{jk}$. Start from the state $(0,0,0)$ at time $0$. At an arrival at node $k$, $k$ carries the trait; at a recovery mark, $k$ loses it; at an arrow from $j$ to $k$, $k$ carries the trait if $j$ does. Since the acquisition rate~\eqref{eq:lambda} is a sum of a term that does not depend on the state and of terms each proportional to one $x_j$, this construction gives the household chain. A \emph{path} from an arrival at $(k_0,s_0)$ to $(j,t)$ is a sequence of times $s_0<s_1<\dots<s_n\le t$ and of nodes $k_0,\dots,k_n=j$ such that at each $s_i$, $i\ge1$, there is an arrow from $k_{i-1}$ to $k_i$, and no recovery mark at $k_i$ lies in $(s_i,s_{i+1}]$, with $s_{n+1}=t$ ($n=0$ is allowed). By induction along the sequence of events, $X_j(t)=1$ if and only if there is a path to $(j,t)$ from some arrival at a time in $(0,t]$; an arrival at a node that already carries the trait does not change the state, but it counts if it has a path. Whether a path exists from a given arrival depends on the arrows and recovery marks only, not on the other arrivals. Let $\mathcal N$ be the set of arrivals with a path to $(j,t)$. Now keep each arrival independently with probability $\vartheta$; the kept arrivals are Poisson processes with rates $\vartheta a_k$, so the construction with them gives the chain with outside rates $\vartheta a$, and its node $j$ carries the trait at time $t$ if and only if some arrival of $\mathcal N$ is kept. Conditionally on the Poisson processes, this has probability $1-(1-\vartheta)^{|\mathcal N|}\ge\vartheta\,\mathbf 1\{\mathcal N\neq\emptyset\}$. Taking expectations, $\Pb(X^{\vartheta}_j(t)=1)\ge\vartheta\,\Pb(X_j(t)=1)$ for every $t$, and letting $t\to\infty$ gives the stationary statement.
\end{proof}

\begin{theorem}[Unique equilibrium]\label{thm:unique-eq}
Under Assumption~\ref{ass:main}, for every constant $u\in[0,\bar u]^3$ the function $s\mapsto g_u(s)/s$ is strictly decreasing on $(0,\infty)$, and $g_u(s)=s$ has exactly one root $s^*(u)\in(0,\hat\omega]$. Hence~\eqref{eq:mf} with the constant intervention $u$ has exactly one equilibrium $\pi_u(s^*(u))$, which attracts every solution. The root $s^*(u)$ is a continuous, nonincreasing function of each $u_k$.
\end{theorem}

\begin{proof}
For a fixed $u$, the outside rates of $\pi_u(s)$ are $\tilde a(s)=c+sd$ with $c_k=(1-u_k)\alpha_k^0>0$ and $d_k=(1-u_k)\eta_k\ge0$, and the transmission rates do not depend on $s$. Let $0<s<s'$ and $\vartheta=s/s'$. Then $\tilde a(s)=c+\vartheta s'd\ge\vartheta\tilde a(s')$, with strict inequality in every component since $c_k>0$. The stationary prevalences are nondecreasing in the outside rates (Lemma~\ref{lem:coupling}, letting $t\to\infty$), so by Lemma~\ref{lem:subhom}, $m(\tilde a(s))\ge m(\vartheta\tilde a(s'))\ge\vartheta m(\tilde a(s'))$, and $g_u(s)/s\ge g_u(s')/s'$. The first inequality is strict for the prevalence of each node. Let $a\le a'$ with $a_k<a'_k$, and couple the two chains as in the proof of Lemma~\ref{lem:coupling}, from the same initial state. For $x\le\tilde x$ with $x_k=\tilde x_k=0$, the acquisition rates satisfy $\lambda_k(x)<\lambda'_k(\tilde x)$, with a difference at least $a'_k-a_k$. On $[t-1,t]$ consider the event: one loss event at $k$ that empties node $k$ in both chains in $[t-1,t-\frac12]$, then one acquisition event at $k$ in $(t-\frac12,t-\frac14]$ whose mark lies between the two acquisition rates, and no other event at node $k$ in $[t-1,t]$. Its probability is at least a constant $C(a'_k-a_k)>0$ that does not depend on $t$, and on it node $k$ carries the trait at time $t$ in the second chain but not in the first. Hence $\Pb(X'_k(t)=1)-\Pb(X_k(t)=1)\ge C(a'_k-a_k)$ for $t\ge1$, and letting $t\to\infty$, $m_k(a')-m_k(a)\ge C(a'_k-a_k)>0$. Applied with $a=\vartheta\tilde a(s')$ and $a'=\tilde a(s)$, which differ in every component by at least $(1-\vartheta)\min_kc_k$, and since $\omega\ge0$ has a positive component, this gives $g_u(s)/s-g_u(s')/s'\ge C'(s'-s)/(ss')>0$. If $s_1<s_2$ were two roots, then $1=g_u(s_1)/s_1>g_u(s_2)/s_2=1$, a contradiction; $s=0$ is not a root because $g_u(0)>0$. Existence of a root in $(0,\hat\omega]$ is Theorem~\ref{thm:scalar}(ii); since $g_u(s)/s$ has a negative derivative, $g_u'(s^*(u))<1$ at the root; and attractivity is Theorem~\ref{thm:scalar}(iii). Finally $g_u(s)$ is continuous in $(s,u)$ and nonincreasing in $u$, and $g_u(s)-s$ changes sign only at $s^*(u)$, so $s^*$ is continuous and nonincreasing in $u$.
\end{proof}

Theorem~\ref{thm:unique-eq} holds for any rates satisfying Assumption~\ref{ass:main}, however strong the transmission inside households or the feedback of the community. It also holds, with the same proof, for several types of households, since $g_u$ is then a positive combination of functions of the same form. For each constant intervention, the long-run prevalences are therefore $m^\infty(u)=\pi_u(s^*(u))S$, computed by a bisection on one variable. The derivative $g_u'$ may exceed $1$ away from the root, but it is smaller than $1$ at the root.

%======================================================================
\section{Optimal control}\label{sec:control}
%======================================================================

For $u\in\U$ let $p_u$ be the solution of~\eqref{eq:mf} with $p(0)=p_0$ and $m_u=p_uS$. We minimise
\begin{equation}\label{eq:J}
J(u)=\int_0^T\Big[w\cdot m_u(t)+\frac12\sum_{k}B_k\,u_k(t)^2\Big]dt
\end{equation}
over $u\in\U$, where $w_k\ge 0$ weights the prevalence of node $k$ and $B_k>0$ is the cost of intervention $k$. The first term counts carriage-time, mainly in people; the second is the usual quadratic cost of effort \citep{lenhart2007}, which grows faster than linearly because the last part of a reduction is the hardest to obtain.

\begin{theorem}[Existence]\label{thm:existence}
Under Assumption~\ref{ass:main}, the problem~\eqref{eq:J} has an optimal intervention $u^*\in\U$.
\end{theorem}

\begin{proof}
Let $u^n$ be a minimising sequence. $\U$ is a bounded, closed and convex subset of $L^2(0,T)^3$, so a subsequence converges weakly in $L^2$ to some $u^*\in\U$. The states $p^n=p_{u^n}$ take values in $\Delta$ and are uniformly Lipschitz in $t$ (their derivatives are bounded by $2\bar q$), so by the Arzelà--Ascoli theorem a further subsequence converges uniformly to some continuous $p^*$. The generator is affine in $u$: $Q(m,u)=Q^{(0)}(m)+\sum_ku_kQ^{(k)}(m)$, with $Q^{(k)}$ continuous in $m$. Hence
\[
p^n(t)=p_0+\int_0^t\Big[p^nQ^{(0)}(p^nS)+\sum_ku^n_k\,p^nQ^{(k)}(p^nS)\Big]ds,
\]
and passing to the limit, using the uniform convergence of $p^nQ^{(k)}(p^nS)$ and the weak convergence of $u^n_k$, shows that $p^*=p_{u^*}$. Finally $\int_0^Tw\cdot m_{u^n}\,dt\to\int_0^Tw\cdot m_{u^*}\,dt$, and $u\mapsto\int_0^T\sum_kB_ku_k^2\,dt$ is convex and continuous on $L^2$, hence weakly lower semicontinuous. So $J(u^*)\le\liminf J(u^n)$, and $u^*$ is optimal.
\end{proof}

The Hamiltonian of the problem is, for $p\in\Delta$, an adjoint (column) vector $\lambda\in\R^8$ and $u\in[0,\bar u]^3$,
\[
\mathcal{H}(p,\lambda,u)=w\cdot(pS)+\frac12\sum_kB_ku_k^2+p\,Q(pS,u)\,\lambda .
\]
For $k\in\{H,A,E\}$ let $D_k$ be the generator in which node $k$ acquires the trait at rate $1$ (from the states with $x_k=0$) and nothing else happens, and let $G_{jk}$ be the generator in which node $k$ acquires the trait at rate $x_j$. Then
\[
Q(m,u)=\sum_k(1-u_k)a_k(m)D_k+\sum_{j,k}\theta_j(u)\beta_{jk}G_{jk}+\sum_k\gamma_kR_k ,
\]
with $R_k$ the generator of loss at node $k$ at rate $1$.

\begin{theorem}[Optimality system]\label{thm:pmp}
Let $u^*$ be optimal, $p^*=p_{u^*}$ and $m^*=p^*S$. There is an absolutely continuous $\lambda:[0,T]\to\R^8$ such that
\begin{equation}\label{eq:adjoint}
-\dot\lambda(t)=Sw+Q\big(m^*(t),u^*(t)\big)\lambda(t)+\kappa(t)\,S\omega,\qquad \lambda(T)=0,
\end{equation}
where $\kappa(t)=\sum_k(1-u_k^*(t))\,\eta_k\,p^*(t)D_k\lambda(t)$, and, for almost every $t$,
\begin{equation}\label{eq:uopt}
u_k^*(t)=\min\Big\{\bar u,\ \max\Big\{0,\ \frac{\Psi_k(t)}{B_k}\Big\}\Big\},
\end{equation}
with the switching functions
\begin{align*}
\Psi_H(t)&=a_H(m^*)\sum_{x:\,x_H=0}p^*(t,x)\,\big[\lambda(t,x+e_H)-\lambda(t,x)\big],\qquad \Psi_A \text{ the same with } A\text{ in place of }H,\\
\Psi_E(t)&=a_E(m^*)\sum_{x:\,x_E=0}p^*(t,x)\,\big[\lambda(t,x+e_E)-\lambda(t,x)\big]\\
&\quad+\sum_{j\in\{H,A\}}\beta_{Ej}\sum_{x:\,x_E=1,\,x_j=0}p^*(t,x)\,\big[\lambda(t,x+e_j)-\lambda(t,x)\big].
\end{align*}
\end{theorem}

\begin{proof}
The problem has a fixed horizon, a free final state and control constraints, and the dynamics and the running cost are smooth in $p$ and continuous in $u$. Pontryagin's maximum principle \citep{pontryagin1962,fleming1975} gives an adjoint with $-\dot\lambda=\nabla_p\mathcal H$, $\lambda(T)=0$, and $u^*(t)$ minimising $\mathcal H(p^*(t),\lambda(t),\cdot)$ over $[0,\bar u]^3$ for almost every $t$ (the problem is normal because the final state is free). Since $m=pS$ and $Q$ depends on $m$ only through $a_k(m)=\alpha_k^0+\eta_k\,\omega\cdot m$,
\[
\nabla_p\big[pQ(pS,u)\lambda\big]=Q(m,u)\lambda+S\,\nabla_m\big[pQ(m,u)\lambda\big]=Q(m,u)\lambda+S\omega\sum_k(1-u_k)\eta_k\,pD_k\lambda,
\]
which gives~\eqref{eq:adjoint}. The Hamiltonian is a strictly convex quadratic function of each $u_k$, separately, with $\partial\mathcal H/\partial u_k=B_ku_k+p\,(\partial Q/\partial u_k)\,\lambda$. Here $\partial Q/\partial u_k=-a_kD_k$ for $k=H,A$ and $\partial Q/\partial u_E=-a_ED_E-\sum_{j\in\{H,A\}}\beta_{Ej}G_{Ej}$, and $(D_k\lambda)(x)=\lambda(x+e_k)-\lambda(x)$ if $x_k=0$ and $0$ otherwise, $(G_{Ej}\lambda)(x)=x_E(1-x_j)[\lambda(x+e_j)-\lambda(x)]$. So $-p(\partial Q/\partial u_k)\lambda=\Psi_k$, and the minimiser over $[0,\bar u]$ is the projection~\eqref{eq:uopt}.
\end{proof}

The adjoint has a simple interpretation: $\lambda(t,x)$ is the marginal cost, from time $t$ to the horizon, of one more household in state $x$, and $\lambda(t,x+e_k)-\lambda(t,x)$ is the cost of node $k$ carrying the trait in addition. The switching function $\Psi_k$ is the expected cost of the acquisitions that intervention $k$ prevents per unit of effort, and~\eqref{eq:uopt} sets the effort where its marginal cost $B_ku_k$ equals this benefit.

\begin{theorem}[Every intervention is used]\label{thm:positive}
Assume $w_k>0$ for every node. Let $(p^*,\lambda,u^*)$ satisfy the optimality system of Theorem~\ref{thm:pmp}. Then for every $t\in[0,T)$ and every $x\in\X$ and node $k$ with $x_k=0$,
\[
\lambda(t,x+e_k)-\lambda(t,x)>0 .
\]
Consequently $\Psi_k(t)>0$ and $u_k^*(t)>0$ for every $k$ and every $t\in(0,T)$, while $u^*(T)=0$; and $u^*$ is continuous.
\end{theorem}

\begin{proof}
Fix the optimal $p^*$ and $u^*$, and let $Q(t)=Q(m^*(t),u^*(t))$. For a bounded measurable function $\nu\ge 0$ on $[0,T]$ let $\Lambda[\nu]$ be the solution of the linear equation $-\dot\lambda=Sw+Q(t)\lambda+\nu(t)S\omega$, $\lambda(T)=0$. By the backward Kolmogorov equation,
\[
\Lambda[\nu](t,x)=\E\Big[\int_t^T\big(w\cdot X_s+\nu(s)\,\omega\cdot X_s\big)\,ds\ \Big|\ X_t=x\Big],
\]
where $X$ is the household chain with generator $Q(s)$. Its acquisition rates are nondecreasing in the state, so the coupling constructed in the proof of Lemma~\ref{lem:coupling}, started at time $t$ from $x\le x'$, gives two copies of this chain with $X_s\le X'_s$ for $s\ge t$; since $w,\omega\ge 0$ and $\nu\ge0$, $\Lambda[\nu](t,\cdot)$ is nondecreasing on $\X$. For the strict inequality, couple the chains from $x$ and $x'=x+e_k$: with positive probability, node $k$ of $X'$ keeps the trait while node $k$ of $X$ does not acquire it on a time interval of positive length after $t$, because all rates are bounded. Hence $\Lambda[\nu](t,x+e_k)-\Lambda[\nu](t,x)\ge w_k\,\E\int_t^T\mathbf 1\{X'_{s,k}=1,X_{s,k}=0\}\,ds>0$ for $t<T$.

The adjoint $\lambda$ of Theorem~\ref{thm:pmp} is the fixed point of $\lambda\mapsto\Lambda[\kappa_\lambda]$, with $\kappa_\lambda(t)=\sum_k(1-u_k^*)\eta_k\,p^*D_k\lambda$. This map is affine, and $\|\Lambda[\kappa_\lambda](t)-\Lambda[\kappa_{\lambda'}](t)\|_\infty\le K\int_t^T\|\lambda(s)-\lambda'(s)\|_\infty\,ds$ for a constant $K$, so it is a contraction for the norm $\sup_te^{-2K(T-t)}\|\lambda(t)\|_\infty$, and the iterates $\lambda^{0}=0$, $\lambda^{n+1}=\Lambda[\kappa_{\lambda^n}]$ converge uniformly to $\lambda$. If $\lambda^n(t,\cdot)$ is nondecreasing for every $t$, then $D_k\lambda^n\ge0$, so $\kappa_{\lambda^n}\ge 0$ and $\lambda^{n+1}$ is nondecreasing. Hence $\lambda(t,\cdot)$ is nondecreasing, $\kappa\ge 0$, and the first part, applied with $\nu=\kappa$, gives the strict inequality.

By Proposition~\ref{prop:wellposed}, $p^*(t,x)>0$ for $t>0$, and $a_k(m^*)\ge\alpha_k^0>0$, so $\Psi_k(t)>0$ for $t\in(0,T)$ (for $\Psi_E$ the second sum is also nonnegative), and~\eqref{eq:uopt} gives $u_k^*(t)>0$. At $t=T$, $\lambda(T)=0$, so $\Psi_k(T)=0$ and $u^*(T)=0$. Finally, $\Psi_k$ depends only on $p^*$, $m^*$ and $\lambda$, which are continuous, so~\eqref{eq:uopt} defines a continuous version of $u^*$.
\end{proof}

Theorem~\ref{thm:positive} says that, with a quadratic cost of effort, it is never optimal to switch an intervention off before the horizon, nor to use only one of them.

\begin{remark}\label{rem:cost}
This conclusion depends on the form of the cost. With the quadratic cost~\eqref{eq:J}, the marginal cost of effort is zero at $u_k=0$, so any positive benefit $\Psi_k$ calls for some effort. With a cost linear in $u_k$, or with a fixed cost of starting an intervention, the optimal $u_k$ is zero whenever $\Psi_k$ stays below a threshold, and some interventions may then not be used at all. The monotonicity of the adjoint, and hence $\Psi_k>0$ for $t\in(0,T)$, does not depend on the cost of the interventions.
\end{remark}

The decrease of the controls to zero at $t=T$ is an effect of the finite horizon: after $T$ nothing is counted.

\begin{theorem}[Uniqueness]\label{thm:unique}
Let Assumption~\ref{ass:main} hold. There is a constant $C$, which depends on $\alpha,\beta,\gamma,\varphi,\omega,w$ and $\bar u$ but not on $T$ or on $B$, such that, if $T\le1$ and $CT<\min_kB_k$, the optimal intervention is unique. Under the same condition, the iteration $u^{n+1}=\mathcal T(u^n)$, where $\mathcal T(u)$ is the right-hand side of~\eqref{eq:uopt} computed with the state and the adjoint of $u$, and the relaxed iteration $u^{n+1}=\frac12\big(u^n+\mathcal T(u^n)\big)$ used in Section~\ref{sec:params}, converge to it from any $u^0\in\U$.
\end{theorem}

\begin{proof}
For $u\in\U$ let $p_u$ be the state, $\lambda_u$ the solution of~\eqref{eq:adjoint} with $(p_u,u)$ in place of $(p^*,u^*)$, and $\mathcal T(u)_k=\min\{\bar u,\max\{0,\Psi_k[u]/B_k\}\}$. By Theorem~\ref{thm:pmp}, every optimal intervention is a fixed point of $\mathcal T$. In what follows, $C$ denotes constants of the kind stated. By Proposition~\ref{prop:wellposed}, $\|p_u-p_v\|_\infty\le CTe^{CT}\|u-v\|_\infty$. The adjoint equation is linear with coefficients bounded by $C$ and source bounded by $C$, so $\|\lambda_u\|_\infty\le CTe^{CT}$; its coefficients are Lipschitz in $(p,u)$, so $\|\lambda_u-\lambda_v\|_\infty\le CTe^{CT}\big(\|p_u-p_v\|_\infty+\|u-v\|_\infty\big)\|\lambda_u\|_\infty+CTe^{CT}\|p_u-p_v\|_\infty$. The switching functions are bilinear in $(p,\lambda)$ with coefficients Lipschitz in $p$, so $|\Psi_k[u]-\Psi_k[v]|\le C\big(\|\lambda_u-\lambda_v\|_\infty+\|\lambda_v\|_\infty\|p_u-p_v\|_\infty\big)$. The projection is $1$-Lipschitz. Collecting the bounds, and using $e^{CT}\le e^C$ for $T\le1$, $\|\mathcal T(u)-\mathcal T(v)\|_\infty\le\vartheta\|u-v\|_\infty$ with $\vartheta=CT/\min_kB_k$ after renaming $C$. Under the condition of the theorem $\vartheta<1$, so $\mathcal T$ is a contraction of $\U$ for the uniform norm; it has at most one fixed point, and the iteration converges by Banach's theorem. The relaxed map $\frac12(u+\mathcal T(u))$ is a contraction with constant $\frac12(1+\vartheta)<1$ and has the same fixed points. Since an optimal intervention exists (Theorem~\ref{thm:existence}), it is this fixed point.
\end{proof}

The constant $C$ can be written explicitly from the proof in terms of $\bar q$, $L$, $\hat\eta\hat\omega$, $\|w\|_1$ and the largest outside and transmission rates, but the resulting value is very conservative and we do not use it. The condition is of the type found in \citet{fister1998} and \citet{lenhart2007}: uniqueness for short horizons or costly interventions. For the horizon used in Section~\ref{sec:numerics} it is not satisfied, but the forward--backward sweep converged to the same solution from all the starting points we tried (Section~\ref{sec:params}).

\begin{corollary}[Optimal value and parameters]\label{cor:value}
Let $J^*=\min_{u\in\U}J(u)$, and regard $\alpha^0$, $\eta$, $\omega$, $\beta$ and $\gamma$ as independent parameters. Then $J^*$ is a nondecreasing function of each of $\alpha_k^0$, $\eta_k$, $\omega_j$ and $\beta_{jk}$, and a nonincreasing function of each $\gamma_k$; the same holds if the initial law is replaced by a larger one for $\preceq$.
\end{corollary}

\begin{proof}
For each fixed $u$, Proposition~\ref{prop:params} orders the prevalences, and $w\ge0$, so $J(u)$ is monotone in each parameter in the stated sense; the minimum over $\U$ of a family of functions that are all nondecreasing in a parameter is nondecreasing in it.
\end{proof}

\emph{The long-run problem.} When the horizon is long, it is natural to compare the optimal control with the best constant intervention for the long-run equilibrium. By Theorem~\ref{thm:unique-eq} and the implicit function theorem, which applies because $g_u(s)$ is a rational function of $(s,u)$ with a positive denominator and $g_u'(s^*(u))<1$, $s^*(u)$ is a continuously differentiable function of $u$. The long-run cost
\[
\bar J(u)=w\cdot m^\infty(u)+\frac12\sum_kB_ku_k^2,\qquad m^\infty(u)=\pi_u(s^*(u))S,
\]
is therefore continuously differentiable and has a minimiser on the compact set $[0,\bar u]^3$. Since $p^\infty(u)=\pi_u(s^*(u))$ by uniqueness, Proposition~\ref{prop:equilibrium} shows that each $m^\infty_k$ is nonincreasing in each $u_j$. We now relate the optimal control on a long horizon to this static problem.

\emph{Turnpike.} Let $(p_s,u_s)$ minimise $\ell(p,u)=w\cdot pS+\frac12\sum_kB_ku_k^2$ over the set of equilibria $\{(p,u)\in\Delta\times[0,\bar u]^3:\ pQ(pS,u)=0\}$, and let $\ell_s$ be the minimum; since the equilibrium is unique for each $u$ (Theorem~\ref{thm:unique-eq}), $u_s$ is the long-run optimum and $\ell_s=\min\bar J$. Let $K\subset\Delta$ be a closed set that contains every state reached from $p_0$ by an intervention in $\U$. When $p_0=\bar p$, one can take $K=\{q\in\Delta:\ p_{\min}(\bar u\mathbf 1)\preceq q\preceq\bar p\}$, where $\bar u\mathbf 1$ is the maximal intervention: for every $u\in\U$, Theorem~\ref{thm:comparison} gives $p_u(t)\preceq\bar p$ (compare with $\tilde u=0$) and $p_u(t)\succeq p_{\bar u\mathbf 1}(t)$; by Proposition~\ref{prop:equilibrium}, $p_{\bar u\mathbf 1}(t)$ is nonincreasing, so $p_{\bar u\mathbf 1}(t)\succeq p^\infty(\bar u\mathbf 1)\succeq p_{\min}(\bar u\mathbf 1)$ by Proposition~\ref{prop:minmax}. This set is closed and convex, since it is defined by finitely many linear inequalities. We write $d(p,u)=\big(\|p-p_s\|_2^2+|u-u_s|_2^2\big)^{1/2}$ and $J_T(u)=\int_0^T\ell(p_u(t),u(t))\,dt$ for the cost on the horizon $T$.

\begin{assumption}[Strict dissipativity]\label{ass:dissip}
There are a function $\Sigma$, continuously differentiable on a neighbourhood of $K$, and a nondecreasing function $\varrho:[0,\infty)\to[0,\infty)$ with $\varrho(0)=0$ and $\varrho(r)>0$ for $r>0$, such that
\[
\ell(p,u)-\ell_s+\nabla\Sigma(p)\cdot pQ(pS,u)\ \ge\ \varrho\big(d(p,u)\big)\qquad\text{for all }(p,u)\in K\times[0,\bar u]^3 .
\]
Moreover, the solution $p_{u_s}$ of~\eqref{eq:mf} with the constant intervention $u_s$ satisfies $C_0=\sup_{T>0}\int_0^T\big(\ell(p_{u_s}(t),u_s)-\ell_s\big)dt<\infty$.
\end{assumption}

The function $\Sigma$ is called a storage function; strict dissipativity is the standard sufficient condition for the turnpike property \citep{grune2016,faulwasser2017}.

\begin{theorem}[Turnpike]\label{thm:turnpike}
Under Assumption~\ref{ass:dissip}, let $V_T$ be the optimal value on $[0,T]$ and $(p_T,u_T)$ an optimal pair, and let $C_1=C_0+\max_K\Sigma-\min_K\Sigma$. Then, for every $T>0$:
\begin{enumerate}
\item[(i)] $T\ell_s-(\max_K\Sigma-\min_K\Sigma)\le V_T\le T\ell_s+C_0$, so that $V_T/T\to\ell_s$;
\item[(ii)] $\int_0^T\varrho\big(d(p_T(t),u_T(t))\big)\,dt\le C_1$; hence, for every $r>0$, the time during which $d(p_T(t),u_T(t))>r$ is at most $C_1/\varrho(r)$, whatever $T$;
\item[(iii)] the constant intervention $u_s$ is nearly optimal: $J_T(u_s)-V_T\le C_1$ for every $T$.
\end{enumerate}
\end{theorem}

\begin{proof}
Let $u\in\U$ and $p$ its state, which stays in $K$. Since $\frac{d}{dt}\Sigma(p(t))=\nabla\Sigma(p)\cdot pQ(pS,u)$,
\begin{align*}
J_T(u)&=T\ell_s+\int_0^T\big[\ell-\ell_s+\nabla\Sigma(p)\cdot pQ(pS,u)\big]dt-\big[\Sigma(p(T))-\Sigma(p(0))\big]\\
&\ge T\ell_s+\int_0^T\varrho(d)\,dt-\Big(\max_K\Sigma-\min_K\Sigma\Big).
\end{align*}
With $u=u_T$ this gives the lower bound in (i) and, combined with $V_T\le J_T(u_s)\le T\ell_s+C_0$, statements (ii) and (iii). The bound on the time spent at distance larger than $r$ follows from $\varrho(d)\ge\varrho(r)$ there.
\end{proof}

Theorem~\ref{thm:turnpike} holds for any initial law $p_0$, with $K$ and $C_0$ chosen for it, and for any horizon, and it needs no condition at the final time, so it covers the free final state of our problem. Its assumption is global, and the natural choice does not satisfy it: with the linear storage $\Sigma(p)=\lambda_s\cdot p$, where $\lambda_s$ is the multiplier of the static problem, the inequality fails at some points of $K\times[0,\bar u]^3$, and a numerical search for a quadratic storage did not produce a certificate. The failure of the linear storage is structural. Let $\lambda_s$ be the multiplier of the static problem, defined below. The optimality condition of the static problem in $p$ reads $Sw+Q(m_s,u_s)\lambda_s+\kappa_sS\omega=\ell_s'\mathbf 1$ for a constant $\ell_s'$, with $\kappa_s=\sum_k(1-u_{s,k})\eta_kp_sD_k\lambda_s$, and a direct computation gives, for $p\in\Delta$,
\[\ell(p,u_s)-\ell_s+\lambda_s\cdot pQ(pS,u_s)=\big(\omega\cdot(p-p_s)S\big)\sum_k(1-u_{s,k})\eta_k\,(p-p_s)D_k\lambda_s ,\]
a product of two linear forms in $p-p_s$. It vanishes whenever $\omega\cdot pS=\omega\cdot p_sS$, and identically when $\varphi=0$. A linear storage therefore never penalises a distance in $p$ along this hyperplane, which meets $K$ in points other than $p_s$. A storage that works must use the contraction of the dynamics; we construct one below (Theorem~\ref{thm:dissip}). The local picture, by contrast, can be described completely. Assume that $p_s(x)>0$ for every $x$ and that $u_s$ lies in the interior of $[0,\bar u]^3$, so that the projection in~\eqref{eq:uopt} is inactive near $u_s$. Let $V$ be an orthonormal basis of $\{y\in\R^8:\sum_xy(x)=0\}$, and let $A$ and $B_u$ be the derivatives of $f(p,u)=pQ(pS,u)$ with respect to $p$ and $u$ at $(p_s,u_s)$, written in the coordinates $\xi=V^{\top}(p-p_s)$; assume that $A$ is invertible. Then there is a multiplier $\lambda_s\in\R^8$ of the static problem, unique up to adding a multiple of $\mathbf 1$, such that $Sw+\nabla_p\big(\lambda_s\cdot f(p,u_s)\big)\big|_{p_s}\in\R\mathbf 1$ and $B_ku_{s,k}=\Psi_k(p_s,\lambda_s)$. Since $f(p,u)\cdot\mathbf 1=0$ for every $p$, the component of the adjoint along $\mathbf 1$, which grows linearly in $T-t$, does not enter the equations; we write $\mu=V^{\top}(\lambda-\lambda_s)$ and $v=u-u_s$. Linearising the optimality system of Theorem~\ref{thm:pmp} at $(p_s,u_s,\lambda_s)$ gives
\[
\dot\xi=A\xi+B_uv,\qquad -\dot\mu=H_{\xi\xi}\xi+H_{\xi v}v+A^{\top}\mu,\qquad H_{vv}v+H_{\xi v}^{\top}\xi+B_u^{\top}\mu=0,
\]
where $H_{\xi\xi},H_{\xi v},H_{vv}$ are the second derivatives of $\ell+\lambda_s\cdot f$ in these coordinates; $H_{vv}=\mathrm{diag}(B_k)$, since $f$ is affine in $u$. Eliminating $v$, $\zeta=(\xi,\mu)$ satisfies $\dot\zeta=\mathbb H\zeta$ with the Hamiltonian matrix
\[
\mathbb H=\begin{pmatrix}A-B_uH_{vv}^{-1}H_{\xi v}^{\top}&-B_uH_{vv}^{-1}B_u^{\top}\\-(H_{\xi\xi}-H_{\xi v}H_{vv}^{-1}H_{\xi v}^{\top})&-(A-B_uH_{vv}^{-1}H_{\xi v}^{\top})^{\top}\end{pmatrix}.
\]
Its eigenvalues are symmetric with respect to the imaginary axis, with multiplicities; if none lies on the axis, its stable and unstable invariant subspaces both have dimension $n=7$. For the nonlinear problem, the final condition $\lambda(T)=0$ gives $\mu(T)=-V^{\top}\lambda_s$, which is not small, and $u(T)=0$ lies on the boundary; the linearisation therefore describes the approach to the static optimum and the departure from it, but not the final layer itself.

\begin{theorem}[Exponential turnpike of the linearised problem]\label{thm:lqturnpike}
Assume that $H_{vv}$ is positive definite, that $\mathbb H$ has no eigenvalue on the imaginary axis, and that its stable invariant subspace $E^s$ and unstable invariant subspace $E^u$ are graphs over $\xi$ and over $\mu$ respectively, that is $E^s\cap\{\xi=0\}=E^u\cap\{\mu=0\}=\{0\}$. Let $\sigma>0$ be smaller than the smallest absolute value of the real parts of the eigenvalues of $\mathbb H$, or equal to it if the eigenvalues with this real part are semisimple. Then there are $T_0$ and $C$ such that, for $T\ge T_0$, the linearised system with $\xi(0)=\xi_0$ and $\mu(T)=\mu_T$ has a unique solution, and
\[
|\xi(t)|+|v(t)|+|\mu(t)|\le C\big(e^{-\sigma t}|\xi_0|+e^{-\sigma(T-t)}|\mu_T|\big)\qquad\text{for all }t\in[0,T].
\]
\end{theorem}

\begin{proof}
Let the columns of $\begin{pmatrix}X_s\\M_s\end{pmatrix}$ and $\begin{pmatrix}X_u\\M_u\end{pmatrix}$ be bases of $E^s$ and $E^u$, with $n\times n$ blocks; the graph conditions say that $X_s$ and $M_u$ are invertible. Since $\R^{2n}=E^s\oplus E^u$, every solution can be written $\zeta(t)=e^{\mathbb Ht}\begin{pmatrix}X_s\\M_s\end{pmatrix}a+e^{\mathbb H(t-T)}\begin{pmatrix}X_u\\M_u\end{pmatrix}b$, and there is $K_1$ with $\|e^{\mathbb Ht}|_{E^s}\|\le K_1e^{-\sigma t}$ and $\|e^{-\mathbb Ht}|_{E^u}\|\le K_1e^{-\sigma t}$ for $t\ge0$. The boundary conditions read $X_sa+R_1(T)b=\xi_0$ and $R_2(T)a+M_ub=\mu_T$, where $R_1(T)$ is the $\xi$-part of $e^{-\mathbb HT}\begin{pmatrix}X_u\\M_u\end{pmatrix}$ and $R_2(T)$ the $\mu$-part of $e^{\mathbb HT}\begin{pmatrix}X_s\\M_s\end{pmatrix}$, both of norm at most $K_1'e^{-\sigma T}$. For $T\ge T_0$ this system is a small perturbation of a block-diagonal invertible one, so it has a unique solution with $|a|\le C(|\xi_0|+e^{-\sigma T}|\mu_T|)$ and $|b|\le C(|\mu_T|+e^{-\sigma T}|\xi_0|)$. Then $|\zeta(t)|\le C(e^{-\sigma t}|a|+e^{-\sigma(T-t)}|b|)$, and $e^{-\sigma t}e^{-\sigma T}\le e^{-\sigma(T-t)}$, $e^{-\sigma(T-t)}e^{-\sigma T}\le e^{-\sigma t}$ give the bound for $\zeta$; $v$ is a linear function of $\zeta$.
\end{proof}

The bound says that the solution leaves the initial state and approaches the static optimum at rate $\sigma$, stays there, and leaves it only in a final layer of the same width, whatever the horizon. In Section~\ref{sec:numerics} we check the assumptions of Theorem~\ref{thm:lqturnpike} for our parameters, and we find that the nonlinear optimal controls approach and leave the static optimum at the rate allowed by the theorem.

\emph{A storage function.} The failure of the linear storage suggests adding to it a function that decreases along the dynamics. Assume that $A$ is stable, and let $P$ be the symmetric positive definite solution of $A^{\top}P+PA=-I$, so that $\mathcal V(\xi)=\frac12\xi^{\top}P\xi$ is a Lyapunov function of the linearised dynamics. For $\varepsilon>0$ let
\[
\Sigma_\varepsilon(p)=\lambda_s\cdot p-\varepsilon\,\mathcal V\big(V^{\top}(p-p_s)\big).
\]
Let $\Pi=VV^{\top}$ be the orthogonal projection on $\{y:\sum_xy(x)=0\}$, $|\cdot|$ the Euclidean norm, $\|\cdot\|$ the spectral norm, and define the constants
\begin{align*}
C_2&=\Big(\sum_k c_{2,k}^2\Big)^{1/2},\quad c_{2,k}=a_k(m_s)\,|\Pi D_k\lambda_s|+\eta_k\,|\Pi S\omega|\max_x|(D_k\lambda_s)(x)|+\mathbf 1_{k=E}\,|\Pi G_E\lambda_s|,\\
C_3&=|\Pi S\omega|\sum_k(1-u_{s,k})\eta_k\,|\Pi D_k\lambda_s|,\qquad
C_4=\|P\|\Big(\sum_k c_{4,k}^2\Big)^{1/2},\quad c_{4,k}=\bar a_k\|D_k\|+\mathbf 1_{k=E}\|G_E\|,
\end{align*}
where $G_E=\sum_{j\in\{H,A\}}\beta_{Ej}G_{Ej}$ and $\bar a_k=\alpha_k\big[(1-\varphi)+\varphi\hat\omega/\omega\cdot\bar m\big]$ is the largest value of $a_k$. Finally let $c_K\ge0$ be such that the contraction inequality
\begin{equation}\label{eq:contraction}
-\nabla\mathcal V(\xi)\cdot V^{\top}\big(pQ(pS,u_s)\big)^{\top}\ \ge\ c_K|\xi|^2,\qquad\xi=V^{\top}(p-p_s),
\end{equation}
holds for all $p\in K$. It holds with $c_K=\frac12-\sqrt2\,\|P\|\,|\Pi S\omega|\,\|D_\eta(u_s)\|$, where $D_\eta(u_s)=\sum_k(1-u_{s,k})\eta_kD_k$, since $|p-p_s|\le\sqrt2$ on $\Delta$; when $\varphi=0$ this gives $c_K=\frac12$.

\begin{theorem}[Strict dissipativity]\label{thm:dissip}
Assume that $p_s(x)>0$ for every $x$, that $u_s$ is interior, that $A$ is stable, and that, for some $\varepsilon>0$,
\[
\varepsilon c_K>C_3\qquad\text{and}\qquad (C_2+\varepsilon C_4)^2<2\Big(\min_kB_k\Big)\big(\varepsilon c_K-C_3\big).
\]
Then Assumption~\ref{ass:dissip} holds with the storage $\Sigma_\varepsilon$ and $\varrho(r)=\kappa r^2$, where $\kappa>0$ is the smallest eigenvalue of $\begin{pmatrix}\varepsilon c_K-C_3&-\frac12(C_2+\varepsilon C_4)\\-\frac12(C_2+\varepsilon C_4)&\frac12\min_kB_k\end{pmatrix}$, and the conclusions of Theorem~\ref{thm:turnpike} hold. When $\varphi=0$ and $c_K=\frac12$, the condition holds for some $\varepsilon$ if and only if $\min_kB_k>4C_2C_4$.
\end{theorem}

\begin{proof}
The constant $C_0$ is finite: the solution with the constant intervention $u_s$ converges to $p_s$ (Theorem~\ref{thm:unique-eq}), and eventually exponentially fast since $A$ is stable. Write $\delta p=p-p_s$, $\xi=V^{\top}\delta p$, so that $|\xi|=|\delta p|$, $\delta u=u-u_s$ and $f(p,u)=pQ(pS,u)$, and split
\[
\ell(p,u)-\ell_s+\nabla\Sigma_\varepsilon(p)\cdot f(p,u)=E_1+E_2+E_3+E_4,
\]
with $E_1=\ell(p,u_s)-\ell_s+\lambda_s\cdot f(p,u_s)$, $E_2=\ell(p,u)-\ell(p,u_s)+\lambda_s\cdot\big(f(p,u)-f(p,u_s)\big)$, $E_3=-\varepsilon\nabla\mathcal V(\xi)\cdot V^{\top}(f(p,u)-f(p,u_s))^{\top}$ and $E_4=-\varepsilon\nabla\mathcal V(\xi)\cdot V^{\top}f(p,u_s)^{\top}$. By the product formula above, $E_1=\big(\delta p\cdot\Pi S\omega\big)\sum_k(1-u_{s,k})\eta_k\,\delta p\cdot\Pi D_k\lambda_s\ge-C_3|\xi|^2$. Since $f$ is affine in $u$, $f(p,u)-f(p,u_s)=-\sum_k\delta u_k\big(a_k(m)pD_k+\mathbf 1_{k=E}\,pG_E\big)$, so $E_2=\frac12\sum_kB_k\delta u_k^2+\sum_k\delta u_k\big(B_ku_{s,k}-\Psi_k(p,\lambda_s)\big)$, where $\Psi_k(p,\lambda_s)$ is the switching function of Theorem~\ref{thm:pmp}. The stationarity of the static problem in $u$ gives $B_ku_{s,k}=\Psi_k(p_s,\lambda_s)$, and
\[
\Psi_k(p,\lambda_s)-\Psi_k(p_s,\lambda_s)=a_k(m_s)\,\delta p\cdot D_k\lambda_s+\eta_k\,(\delta p\cdot S\omega)\,(p\cdot D_k\lambda_s)+\mathbf 1_{k=E}\,\delta p\cdot G_E\lambda_s,
\]
whose absolute value is at most $c_{2,k}|\xi|$; hence $E_2\ge\frac12\min_kB_k|\delta u|^2-C_2|\xi||\delta u|$. Next, $|E_3|\le\varepsilon\|P\||\xi|\sum_k|\delta u_k|\,c_{4,k}\le\varepsilon C_4|\xi||\delta u|$, since $|pD_k|\le\|D_k\|$ and $|pG_E|\le\|G_E\|$ for $p\in\Delta$. Finally $E_4\ge\varepsilon c_K|\xi|^2$ by~\eqref{eq:contraction}. Adding,
\[
E_1+E_2+E_3+E_4\ \ge\ (\varepsilon c_K-C_3)|\xi|^2-(C_2+\varepsilon C_4)|\xi||\delta u|+\tfrac12\min_kB_k|\delta u|^2\ \ge\ \kappa\big(|\xi|^2+|\delta u|^2\big),
\]
and the matrix is positive definite exactly under the stated condition. For the bound on $c_K$: the reduced dynamics with $u=u_s$ are $\dot\xi=A\xi+N(\xi)$, where $N(\xi)=V^{\top}\big((\delta p\cdot S\omega)\,\delta p\,D_\eta(u_s)\big)^{\top}$ collects the quadratic terms, because $Q$ is affine in $m$; then $-\xi^{\top}P(A\xi+N(\xi))=\frac12|\xi|^2-\xi^{\top}PN(\xi)$ and $|N(\xi)|\le|\Pi S\omega|\,\|D_\eta(u_s)\|\,|\xi|^2$, so that $|\xi^{\top}PN(\xi)|\le\sqrt2\,\|P\|\,|\Pi S\omega|\,\|D_\eta(u_s)\|\,|\xi|^2$. When $\varphi=0$, $C_3=0$; with $c_K=\frac12$ the condition $(C_2+\varepsilon C_4)^2<\min_kB_k\,\varepsilon$ holds for some $\varepsilon$ if and only if $\min_kB_k>\min_{\varepsilon>0}(C_2+\varepsilon C_4)^2/\varepsilon=4C_2C_4$.
\end{proof}

Theorem~\ref{thm:dissip} gives the turnpike property, without restriction on the horizon or on the final state, when the interventions are costly enough compared with the sensitivity of the switching functions ($C_2$, $C_4$) and when the community feedback ($C_3$, and its effect on $c_K$) is weak enough. The explicit constants are conservative. For a given set of parameters, the inequality of Assumption~\ref{ass:dissip} with the storage $\Sigma_\varepsilon$ is a polynomial inequality on the polytope $K\times[0,\bar u]^3$, and it can be checked by a sum-of-squares certificate \citep{parrilo2003,lasserre2001}, which we do in Section~\ref{sec:numerics}.

\emph{Uncertain parameters.} Some parameters, such as $\varphi$ and $c$, are not known. Let $\theta_1,\dots,\theta_n$ be scenarios for the uncertain parameters, with weights $\varpi_i>0$ summing to $1$, for example the nodes and weights of a quadrature rule for a prior distribution. A single intervention $u$ is applied, whatever the scenario, and we minimise the expected cost
\[
J_{\varpi}(u)=\sum_{i=1}^n\varpi_i\,J_{\theta_i}(u),
\]
where $J_{\theta_i}$ is the cost~\eqref{eq:J} computed with the parameters $\theta_i$ and the initial law of scenario $i$. This is the problem of Section~\ref{sec:control} for the product system made of the $n$ mean-field equations, which share the control but not the state. Theorems~\ref{thm:existence} to~\ref{thm:unique} therefore apply to it with the same proofs: an optimal intervention exists; it satisfies $u_k^*=\min\{\bar u,\max\{0,\sum_i\varpi_i\Psi^{i}_k/B_k\}\}$, where $\Psi^i_k$ is the switching function of scenario $i$, computed with its own state and adjoint; it is positive on $(0,T)$ when $w_k>0$, since each $\Psi^i_k$ is positive; and it is unique for short horizons. The \emph{expected value of perfect information},
\[
\mathrm{EVPI}=\min_{u}J_\varpi(u)-\sum_i\varpi_i\min_{u}J_{\theta_i}(u)\ \ge0,
\]
is the expected gain that knowing the true scenario before choosing the intervention would bring; it measures what a study that estimates the uncertain parameters is worth for the choice of the intervention.

\emph{Several types of households.} The households of a town differ in their exposure. Let there be $Z$ types, in proportions $\rho_1,\dots,\rho_Z$, with outside rates $\alpha^{(z)}$ such that $\alpha^{(z)}_k>0$ for every $z$ and $k$, and let all types share the community, so that $m=\sum_z\rho_zm^{(z)}$ and the community term of type $z$ is $a^{(z)}_k(m)=\alpha^{(z)}_k\big[(1-\varphi)+\varphi\,\omega\cdot m/\omega\cdot\bar m\big]=\alpha^{(z),0}_k+\eta^{(z)}_k\,\omega\cdot m$, with $\eta^{(z)}_k=\varphi\alpha^{(z)}_k/(\omega\cdot\bar m)$. The mean-field equation becomes a system for the laws $p^{(z)}$ of the $Z$ types, coupled through $m$. All the results of Sections~\ref{sec:mf}--\ref{sec:control} remain true, with the same proofs, when the state of a household is the pair $(z,x)$ and the order compares states of the same type; the constants change, for example $\sqrt8$ becomes $\sqrt{8Z}$ in Theorem~\ref{thm:mf}, and the constant of Theorem~\ref{thm:unique} also depends on $\min_z\rho_z$; in Theorem~\ref{thm:scalar}, $g_u(s)=\sum_z\rho_z\,\omega\cdot\pi^{(z)}_u(s)S$. The interventions may be \emph{targeted}, with a different $u^{(z)}$ for each type, or \emph{uniform}. With the cost $\frac12\sum_z\rho_z\sum_kB_k(u_k^{(z)})^2$, which is proportional to the number of households covered, the optimality conditions become
\[
u_k^{(z)}=\min\Big\{\bar u,\max\Big\{0,\frac{\Psi_k^{(z)}}{\rho_zB_k}\Big\}\Big\}\ \text{(targeted)},\qquad u_k=\min\Big\{\bar u,\max\Big\{0,\frac{\sum_z\Psi_k^{(z)}}{B_k}\Big\}\Big\}\ \text{(uniform)},
\]
where $\Psi^{(z)}_k$ is the switching function of Theorem~\ref{thm:pmp} computed with the mass $\rho_zp^{(z)}$ of type $z$ and its adjoint $\lambda^{(z)}$, which satisfies~\eqref{eq:adjoint} with $Q$ replaced by the generator of type $z$ and $\kappa=\sum_{z'}\sum_k(1-u_k^{(z')})\eta_k^{(z')}\rho_{z'}p^{(z')}D_k\lambda^{(z')}$.

%======================================================================
\section{Parameters and numerical method}\label{sec:params}
%======================================================================

The parameters (Table~\ref{tab:params}) are those of a typical resistance trait in Tenkodogo, estimated from a survey of $59$ households in which a person, an animal and the soil of the yard were sampled (Appendix~\ref{app:estimation}). Time is measured in units of the mean duration of carriage, taken equal for the three nodes; for resistant \emph{E.~coli} in people it is about $110$ days in household members in the Netherlands \citep{haverkate2017} and about $102$ days in rural Burkina Faso \citep{aizouk2025}; we use $110$ days to convert time to months. The rates of acquisition from outside of people and animals are $\alpha_H=0.31$ and $\alpha_A=0.34$, and the coupling between person and animal is $c=0.53$. A survey at one time gives the joint law of the two nodes, which is symmetric in them; it cannot tell whether transmission goes mostly from the animal to the person or the reverse, so we use the same rate in both directions ($\beta_{HA}=\beta_{AH}=c$). Soil isolates were available in few households, and the survey cannot separate the acquisition of the soil from outside from transmission to it. We take $\alpha_E=0.15$, a rate $\beta_E=0.3$ from the soil to people and animals and a rate $\beta_{\to E}=0.5$ from people and animals to the soil, and we vary all three. With these values the baseline soil prevalence is $0.30$, below the $0.44$ observed in the $70$ soil isolates; a scenario calibrated on the observed value is given in Section~\ref{sec:sens}. The share of acquisition due to the community, $\varphi$, cannot be estimated from one survey either. In a source attribution study, human-to-human transmission within or between households accounted for about $60\%$ of community-acquired carriage of resistant \emph{E.~coli}, and food, animals and the environment for the rest \citep{mughinigras2019}; part of these sources depends on the prevalence in the community and part does not. We take $\varphi=0.5$ as central value, $0.25$ to $0.75$ as plausible range, and vary it from $0$ to $0.9$. The weights $\omega$ are equal.

The cost weights express the priority of public health: $w=(1,0.5,0.25)$, so that human carriage counts twice as much as animal carriage and four times as much as soil contamination. The costs of the interventions are not known in money; we take $B_H=B_A=B_E=1$, which means that a reduction of $50\%$ of a route kept for one time unit costs as much as $0.125$ units of human carriage-time, and we study a range of costs from $0.2$ to $5$ times this value. The horizon is $T=10$ time units, about three years, and $\bar u=0.9$. The initial state is the baseline equilibrium $\bar p$; with these parameters the baseline prevalences are $\bar m=(0.35,\,0.36,\,0.30)$.

\emph{Survey sample and mean-field model.} The rates are estimated from $59$ households, while the mean-field equation describes a community of many households. There is no contradiction: the $59$ households are a sample from a town of several thousand households, the rates are properties of the household chain, which the sample estimates, and the mean-field equation is applied to the town, not to the sample. Theorem~\ref{thm:mf} and the simulations of Section~\ref{sec:numerics} show that the mean-field equation describes a community of a few thousand households well. Its use for a village of a hundred households would be questionable, and the estimates themselves are uncertain because the sample is small; we propagate this uncertainty to the optimal strategy by a bootstrap in Section~\ref{sec:sens}.

\begin{table}[t]
\centering
\caption{Parameters. Time unit: mean duration of human carriage (about $110$ days). Source: S, estimated from the survey (Appendix~\ref{app:estimation}); L, literature; A, assumption (varied in the sensitivity study); P, policy choice.}
\label{tab:params}
\small
\begin{tabular}{p{3.3cm}p{6.2cm}p{3.6cm}l}
\toprule
Parameter & Meaning & Value & Source\\
\midrule
$\alpha_H,\alpha_A$ & acquisition from outside at baseline & $0.31,\ 0.34$ & S\\
$\alpha_E$ & acquisition of the soil from outside & $0.15$ & A\\
$c=\beta_{HA}=\beta_{AH}$ & person--animal transmission & $0.53$ & S\\
$\beta_E=\beta_{EH}=\beta_{EA}$ & soil to people and animals & $0.3$ & A\\
$\beta_{\to E}=\beta_{HE}=\beta_{AE}$ & people and animals to soil & $0.5$ & A\\
$\gamma_H,\gamma_A,\gamma_E$ & loss of carriage (time unit about $110$ days) & $1,\ 1,\ 1$ & L, A\\
$\varphi$ & share of outside acquisition due to the community & $0.5$ & A\\
$\omega$ & weights of community prevalences & $(1,1,1)$ & A\\
$w$ & weights of prevalences in the cost & $(1,\ 0.5,\ 0.25)$ & P\\
$B_H,B_A,B_E$ & costs of the interventions & $1,\ 1,\ 1$ & P\\
$\bar u$, $T$ & bound on the interventions, horizon & $0.9$, $10$ (about $3$ years) & P\\
\bottomrule
\end{tabular}
\end{table}

\emph{Numerical method.} The state equation~\eqref{eq:mf} and the adjoint equation~\eqref{eq:adjoint} were solved with the classical fourth-order Runge--Kutta scheme on $400$ steps, forward and backward, and the optimality system with the forward--backward sweep \citep{lenhart2007}: given $u^n$, compute $p^n$, then $\lambda^n$, then $\mathcal T(u^n)$ by~\eqref{eq:uopt}, and set $u^{n+1}=\frac12(u^n+\mathcal T(u^n))$, until $\sup|\mathcal T(u^n)-u^n|<10^{-7}$. For the baseline this took $26$ iterations. We checked the method in three ways. The derivative of $J$ given by the adjoint, $B_ku_k-\Psi_k$, agreed with finite differences of $J$ to a relative error of $10^{-5}$ for perturbations of each control on a time window. Refining the grid to $1\,600$ steps changed the optimal cost by $8\times10^{-6}$. And the sweep converged to the same control, within $10^{-7}$, from three different starting points (no intervention, maximal intervention, random intervention).

For the other computations of Section~\ref{sec:numerics}: the function $g_u$ of Theorem~\ref{thm:scalar} was evaluated by solving the $8\times8$ stationary equations, and its roots were located by sign changes on a grid of $401$ points of $[0,\hat\omega]$ ($201$ points for the parameter sets of the sensitivity study and $301$ for the search for several roots) and refined by bisection; its slope was computed by finite differences on the same grid. The long-run optimum was found by the Nelder--Mead method over the three constant interventions. The optimal control of the compartmental model was computed by minimising its discretised cost over controls constant on $40$ intervals with the L-BFGS-B method. The two-zone model was solved by the same forward--backward sweep, with $29$ iterations for the targeted and $26$ for the uniform strategy.

%======================================================================
\section{Numerical results}\label{sec:numerics}
%======================================================================

\subsection{The optimal strategy}

\begin{figure}[t]
\centering
\includegraphics[width=\textwidth]{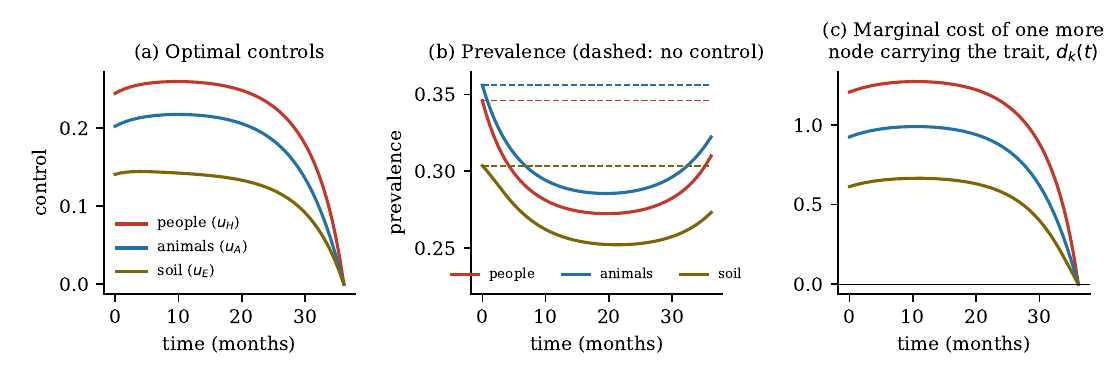}
\caption{Optimal strategy for the baseline parameters. (a) Optimal interventions. (b) Prevalences under the optimal interventions (solid) and without intervention (dashed). (c) Marginal cost of one more node carrying the trait, $d_k(t)=\sum_{x_k=0}p^*(t,x)[\lambda(t,x+e_k)-\lambda(t,x)]/\sum_{x_k=0}p^*(t,x)$.}
\label{fig:baseline}
\end{figure}

Figure~\ref{fig:baseline} shows the optimal interventions and their effect. All three interventions are used from the start, as Theorem~\ref{thm:positive} predicts. The intervention on people is the strongest ($u_H^*$ about $0.25$, a reduction by a quarter of the acquisition of people from outside), followed by the intervention on animals ($u_A^*\approx0.21$) and on the soil ($u_E^*\approx0.14$). The controls stay close to their maximum for about a year and a half and then decrease to zero at the horizon. Human prevalence falls from $0.35$ to a minimum of $0.27$ after about $20$ months; over the three years, $17\%$ of human carriage-time is averted, $16\%$ in animals and $13\%$ in the soil, and the total cost $J$ is reduced by $7.7\%$. Panel (c) shows the marginal cost of one more node carrying the trait, averaged over the households in which that node is free: it is positive at all times, as Theorem~\ref{thm:positive} states, and largest for people, because human carriage is weighted most and because people feed the other nodes.

\subsection{Comparison of strategies}

\begin{figure}[t]
\centering
\includegraphics[width=\textwidth]{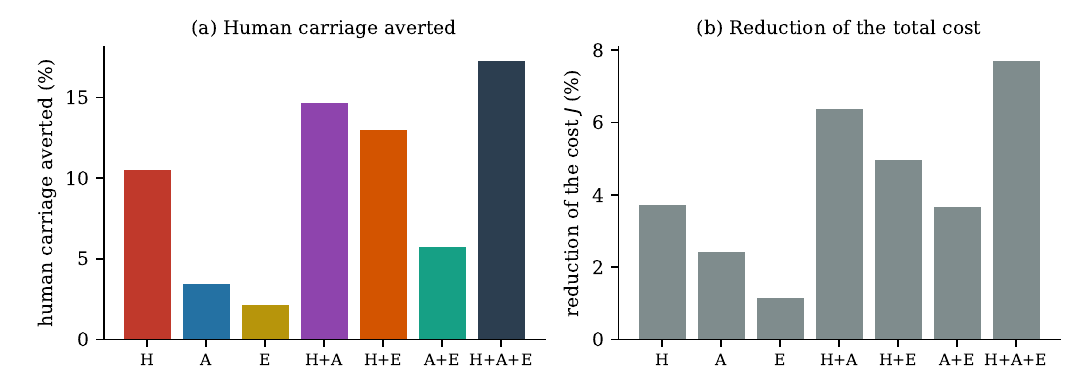}
\caption{Comparison of strategies: (a) human carriage-time averted over three years and (b) reduction of the total cost $J$, for the optimal strategy with each set of interventions.}
\label{fig:strategies}
\end{figure}

We solved the optimal control problem with each intervention alone, with each pair, and with all three (Table~\ref{tab:strategies} and Figure~\ref{fig:strategies}). The intervention on people alone averts $10.5\%$ of human carriage, the intervention on animals alone $3.4\%$ and the intervention on the soil alone $2.1\%$. Combined, the three interventions avert $17.3\%$, more than the sum of the three single interventions ($16.1\%$); the pair on people and animals averts $14.7\%$, more than $10.5+3.4=13.9\%$. Since the single strategies use different efforts, we also applied each component of the optimal combined intervention alone: the three components avert $11.0\%$, $3.6\%$ and $2.3\%$, together $16.9\%$, again less than the $17.3\%$ of the combination. The synergy comes from the transmission inside the household. Without it ($c=\beta_E=\beta_{\to E}=0$), the combination averts less than the sum of its components applied alone ($23.3\%$ against $23.8\%$), because the community feedback alone gives diminishing returns; without community feedback ($\varphi=0$) but with household transmission, it averts more ($10.1\%$ against $9.9\%$). An intervention on animals lowers the transmission from animals to people inside the household, so that the intervention on people prevents acquisitions that would otherwise have been passed on to the animals and back. Over $60$ of the parameter sets of the global sensitivity study, the combined strategy averted more than the sum of the three single strategies in $75\%$ of the cases, so the synergy depends on the parameters. Animal and soil interventions matter mostly for animal and soil carriage, which they reduce by $8.8\%$ and $3.7\%$ on their own. Finally, the optimal time-varying strategy is only slightly better than the best constant strategy, $(u_H,u_A,u_E)=(0.22,0.18,0.12)$ over the three years, which averts $15.9\%$ of human carriage and reduces $J$ by $7.1\%$, instead of $17.3\%$ and $7.7\%$.

\begin{table}[t]
\centering
\caption{Optimal strategies with one, two or three interventions. Carriage-time averted over three years, in people, animals and soil; reduction of the total cost $J$; and mean intervention over the horizon.}
\label{tab:strategies}
\small
\begin{tabular}{lccccc}
\toprule
Interventions & \multicolumn{3}{c}{Carriage averted (\%)} & $J$ reduced & Mean $(u_H,u_A,u_E)$\\
 & people & animals & soil & (\%) & \\
\midrule
H & 10.5 & 3.8 & 4.5 & 3.7 & (0.21, 0.00, 0.00)\\
A & 3.4 & 8.8 & 3.9 & 2.4 & (0.00, 0.17, 0.00)\\
E & 2.1 & 2.0 & 3.7 & 1.2 & (0.00, 0.00, 0.11)\\
H+A & 14.7 & 13.4 & 9.0 & 6.4 & (0.22, 0.18, 0.00)\\
H+E & 13.0 & 6.1 & 8.5 & 5.0 & (0.22, 0.00, 0.12)\\
A+E & 5.7 & 11.1 & 7.9 & 3.7 & (0.00, 0.17, 0.12)\\
H+A+E & 17.3 & 15.9 & 13.3 & 7.7 & (0.22, 0.18, 0.12)\\
\bottomrule
\end{tabular}
\end{table}

\subsection{Cost--effectiveness frontier}

\begin{figure}[t]
\centering
\includegraphics[width=0.55\textwidth]{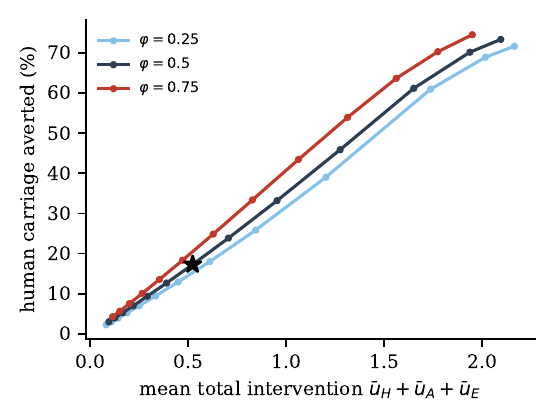}
\caption{Cost--effectiveness frontier: human carriage-time averted over three years by the optimal strategy against the mean total intervention, for costs scaled from $0.2$ to $5$ times the baseline, and three values of the community share $\varphi$. Star: baseline.}
\label{fig:frontier}
\end{figure}

The costs $B_k$ are relative and their units are arbitrary. To remove this dependence, we scaled all costs by a common factor from $0.2$ to $5$ and plotted the human carriage averted against the total effort, measured by the mean total intervention $\bar u_H+\bar u_A+\bar u_E$ over the horizon, which does not depend on the units of $B$ (Figure~\ref{fig:frontier}). Each point of a curve is an optimal strategy for some level of the costs, and a health service can read on the curve what a given effort buys. With $\varphi=0.5$, a total effort of $0.5$ (as in the baseline strategy, with mean reductions of $0.22$, $0.18$ and $0.12$ of the three routes) averts about $17\%$ of human carriage, and a total effort of $1.3$ about $46\%$. The benefit grows slightly faster than the effort, because of the community feedback, and it depends on $\varphi$: for the same total effort of about $0.5$ it is about $15\%$ with $\varphi=0.25$ and $20\%$ with $\varphi=0.75$.

\subsection{Mean-field limit}

\begin{figure}[t]
\centering
\includegraphics[width=\textwidth]{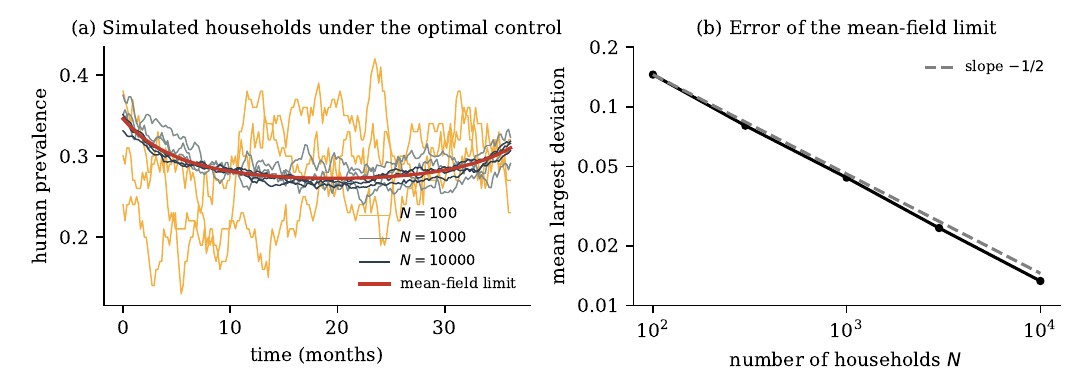}
\caption{Mean-field limit. (a) Human prevalence in three simulated communities for each of $N=100$, $1\,000$ and $10\,000$ households under the optimal intervention, and the mean-field solution. (b) Mean over replicates of the largest deviation $\sup_t\max_k|m^N_k(t)-m_k(t)|$, against $N$; the dashed line has slope $-1/2$.}
\label{fig:mf}
\end{figure}

To check Theorem~\ref{thm:mf} at realistic community sizes, we simulated exactly, with Gillespie's algorithm, communities of $N$ households under the optimal intervention of Figure~\ref{fig:baseline}, starting from independent households with law $\bar p$. Because households are exchangeable, the simulation follows the counts of households in each of the eight states. The deviations were computed on a grid of $201$ times. Figure~\ref{fig:mf} shows that the simulated human prevalence follows the mean-field solution closely from a few thousand households, and that the error decreases as $N^{-1/2}$, as in Theorem~\ref{thm:mf}: the mean of the largest deviation of the three prevalences over three years was 0.145 with $100$ households, 0.044 with $1\,000$ and 0.013 with $10\,000$. A town of several thousand households is well described by the mean-field equation; a village of $100$ households is not, and its prevalences can deviate by more than $0.1$ from the mean-field values, which is as large as the effect of the interventions.

\subsection{Sensitivity study}\label{sec:sens}

\begin{figure}[t]
\centering
\includegraphics[width=\textwidth]{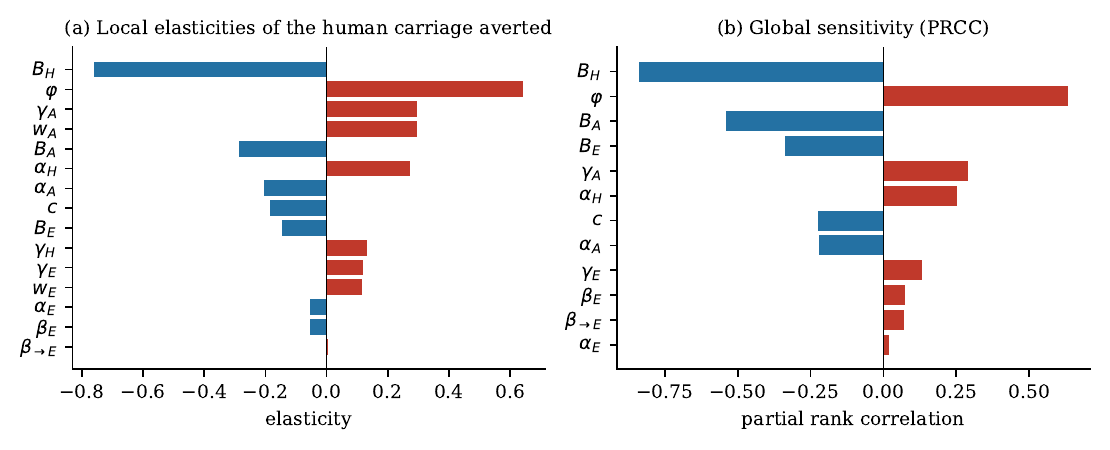}
\caption{Sensitivity of the human carriage averted by the optimal strategy. (a) Local elasticities at the baseline parameters. (b) Partial rank correlation coefficients over $250$ parameter sets drawn by Latin hypercube sampling.}
\label{fig:sens}
\end{figure}

Figure~\ref{fig:sens}a gives the local elasticities of the human carriage averted by the optimal strategy, that is the relative change of the output divided by the relative change of each parameter, computed by central differences with changes of $\pm10\%$. The output is most sensitive to the cost of the human intervention ($-0.76$) and to the share of acquisition due to the community ($0.64$), then to the duration of carriage in animals ($\gamma_A$, $0.30$), to the cost of the animal intervention ($-0.29$), to the weight of animal carriage ($0.29$), to the rate of acquisition of people from outside ($0.27$) and, negatively, to that of animals ($-0.21$). The transmission and acquisition rates of the soil, which are the least known, matter least ($|{\rm elasticity}|\le0.12$; $-0.15$ for the cost of the soil intervention). A stronger coupling between person and animal slightly reduces the proportion of human carriage averted ($-0.18$), because more human carriage then comes from inside the household, where the interventions act only indirectly.

For a global analysis, we drew $250$ parameter sets by Latin hypercube sampling \citep{mckay1979} over wide ranges: the rates $\alpha_k$ from half to twice their values, $c$ from $0$ to $1.5$, $\beta_E$ and $\beta_{\to E}$ from $0$ to $1$, $\gamma_A$ and $\gamma_E$ from $0.5$ to $4$, $\varphi$ from $0$ to $0.9$ and the costs $B_k$ from $0.25$ to $4$, all on a log scale except $c$, $\beta_E$, $\beta_{\to E}$ and $\varphi$. We solved the optimal control problem for each set and computed partial rank correlation coefficients \citep{marino2008}. Over these ranges, the human carriage averted by the optimal strategy had median $22\%$ (5th and 95th percentiles $6\%$ and $66\%$). The partial rank correlations (Figure~\ref{fig:sens}b) confirm the local analysis: the cost of the human intervention ($-0.84$) and the community share ($0.63$) come first, followed by the costs of the animal and soil interventions ($-0.54$ and $-0.34$), the clearance rate in animals ($0.29$), the acquisition rate of people ($0.25$), the coupling ($-0.23$) and the acquisition rate of animals ($-0.22$); the transmission rates involving the soil have partial rank correlations below $0.1$. Each optimal intervention is driven mainly by its own cost, by the parameters of its own route and by the community share: $u_H^*$ by $B_H$ ($-0.99$), $\alpha_H$ ($0.89$) and $\varphi$ ($0.54$); $u_A^*$ by $B_A$ ($-0.98$), $\alpha_A$ ($0.80$), the clearance rate $\gamma_A$ ($-0.76$), $\varphi$ ($0.62$) and the coupling $c$ ($0.53$); $u_E^*$ by $B_E$ ($-0.95$), $\beta_E$ ($0.83$), $\gamma_E$ ($-0.79$), $\alpha_E$ ($0.68$) and $\varphi$ ($0.55$). The order of the three interventions is therefore not fixed: the human intervention was stronger than the animal intervention in $63\%$ of the parameter sets and the strongest of the three in $48\%$, and the order people, animals, soil of the baseline held in only $27\%$. Finally, the weight $\omega_E$ of the soil in the community pressure matters little: with $\omega_E=0$, $0.5$, $1$ and $2$ (and $\omega_H=\omega_A=1$), the human carriage averted is $18.1\%$, $17.6\%$, $17.3\%$ and $16.9\%$.

\emph{Soil calibrated on the survey.} With the outside rate of the soil chosen so that the baseline soil prevalence equals the $0.44$ observed in the survey ($\alpha_E=0.49$, other parameters unchanged), the baseline prevalences become $(0.37,\,0.38,\,0.44)$, and the optimal strategy averts $17.0\%$ of human carriage instead of $17.3\%$. The effort moves towards the soil: the mean interventions are $(0.20,\,0.15,\,0.21)$ instead of $(0.22,\,0.18,\,0.12)$. The main results do not depend on this choice, but the effort on the soil does.

\begin{figure}[t]
\centering
\includegraphics[width=\textwidth]{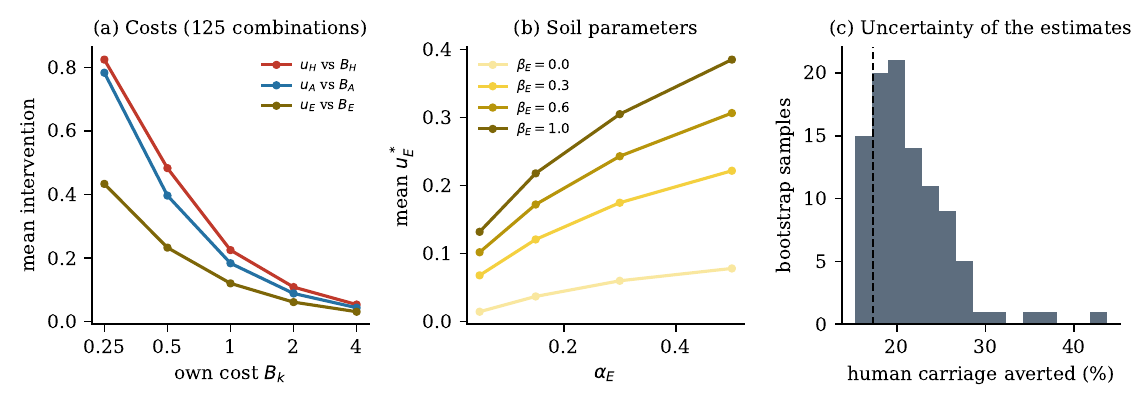}
\caption{(a) Mean optimal intervention against its own cost, averaged over the $125$ combinations of $(B_H,B_A,B_E)$. (b) Mean optimal intervention on the soil against $\alpha_E$, for four values of $\beta_E$ (averaged over $\beta_{\to E}$). (c) Human carriage averted by the optimal strategy for $100$ bootstrap samples of the survey; dashed line: baseline.}
\label{fig:review}
\end{figure}

\emph{Costs of the three interventions.} The costs $B_k$ are the least known inputs of the policy model. We solved the optimal control problem on a grid of $125$ combinations, each $B_k$ taking the values $0.25$, $0.5$, $1$, $2$ and $4$ (Figure~\ref{fig:review}a). The human carriage averted ranges from $4\%$ to $71\%$. Each optimal intervention depends mainly on its own cost: averaged over the other two costs, $u_H^*$ goes from $0.82$ to $0.05$ as $B_H$ goes from $0.25$ to $4$, and $u_A^*$ and $u_E^*$ behave in the same way. The human intervention is the strongest in $47\%$ of the combinations, and the baseline order (people, animals, soil) holds in $28\%$. The allocation of effort therefore cannot be recommended without local costs; what the model gives is the rule that each intervention should be set from its own cost and from the acquisition rate of its route.

\emph{Soil parameters.} The rates involving the soil are assumed. We solved the problem on a grid of $48$ combinations of $\alpha_E\in\{0.05,0.15,0.3,0.5\}$, $\beta_E\in\{0,0.3,0.6,1\}$ and $\beta_{\to E}\in\{0,0.5,1\}$ (Figure~\ref{fig:review}b). The human carriage averted stays between $17\%$ and $25\%$ and the reduction of $J$ between $7\%$ and $12\%$, but the intervention on the soil varies widely, from $0.01$ to $0.39$. It grows with $\alpha_E$ (mean $0.08$ to $0.25$) and with $\beta_E$ (mean $0.05$ to $0.26$), and hardly depends on $\beta_{\to E}$. The conclusions on human carriage are thus robust to the soil parameters, but the recommended effort on the soil is not, and it requires data on the soil.

\emph{Uncertainty of the estimates.} We resampled the $59$ households $100$ times with replacement, re-estimated $\alpha_H$, $\alpha_A$ and $c$ as in Appendix~\ref{app:estimation}, and recomputed the optimal strategy each time. The $95\%$ bootstrap intervals are wide: $0.01$ to $0.59$ for $\alpha_H$, $0.01$ to $0.75$ for $\alpha_A$ and $0$ to $1.25$ for $c$. The human carriage averted by the optimal strategy lies between $16\%$ and $33\%$ (Figure~\ref{fig:review}c), and the optimal human and animal interventions between $0.01$ and $0.41$ and between $0.01$ and $0.35$. With $59$ households, the survey fixes the order of magnitude of the benefit better than the allocation of effort between people and animals.

\subsection{One-dimensional studies}

\begin{figure}[t]
\centering
\includegraphics[width=\textwidth]{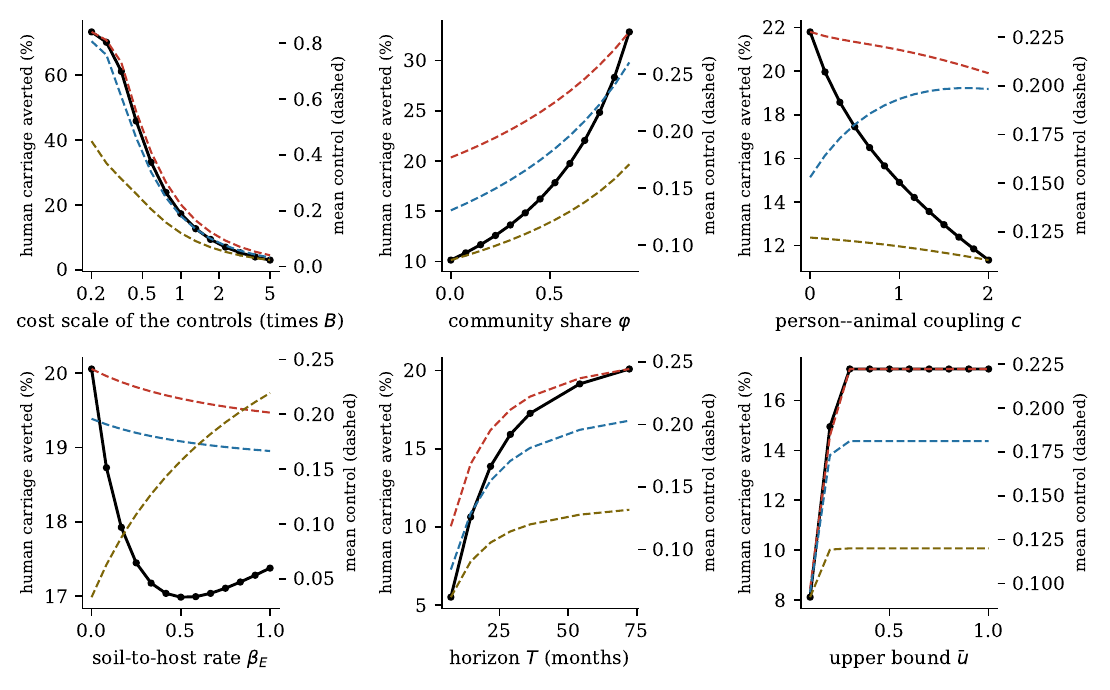}
\caption{One-dimensional studies. Human carriage-time averted by the optimal strategy (black, left axis) and mean optimal interventions (dashed, right axis; red people, blue animals, brown soil) as functions of the scale of the costs $B$, the community share $\varphi$, the coupling $c$, the rate $\beta_E$ from the soil, the horizon $T$ and the bound $\bar u$.}
\label{fig:oned}
\end{figure}

Figure~\ref{fig:oned} shows the effect of six quantities, one at a time. The cost of the interventions has the largest effect: when all costs are divided by five, the human carriage averted rises from $17\%$ to $73\%$, and the interventions on people and animals reach their upper bound $\bar u$ in the middle of the horizon; when they are multiplied by five, it falls to $3\%$. The share $\varphi$ of acquisition due to the community raises both the optimal effort and its effect: with $\varphi=0$, that is without community feedback, $10\%$ of human carriage is averted, and with $\varphi=0.9$, $33\%$. A model that ignores the community therefore underestimates the benefit of the interventions. A stronger coupling $c$ shifts the effort from people towards animals and lowers the proportion averted. The rate $\beta_E$ from the soil changes mainly the intervention on the soil: without transmission from the soil, $u_E^*$ is small, and as $\beta_E$ grows $u_E^*$ grows while the other two decrease slightly. The proportion of human carriage averted is not monotone in $\beta_E$: it is largest without transmission from the soil and smallest near $\beta_E=0.5$, because $\beta_E$ raises both the baseline carriage and the reach of the soil intervention. A longer horizon gives a larger proportion averted (from $5.5\%$ over seven months to $20\%$ over six years), because the slow decrease of the reservoir has a larger effect over a longer horizon. The bound $\bar u$ matters only below $0.3$, since the unconstrained optimum stays below it.

\subsection{Two-dimensional studies}

\begin{figure}[t]
\centering
\includegraphics[width=\textwidth]{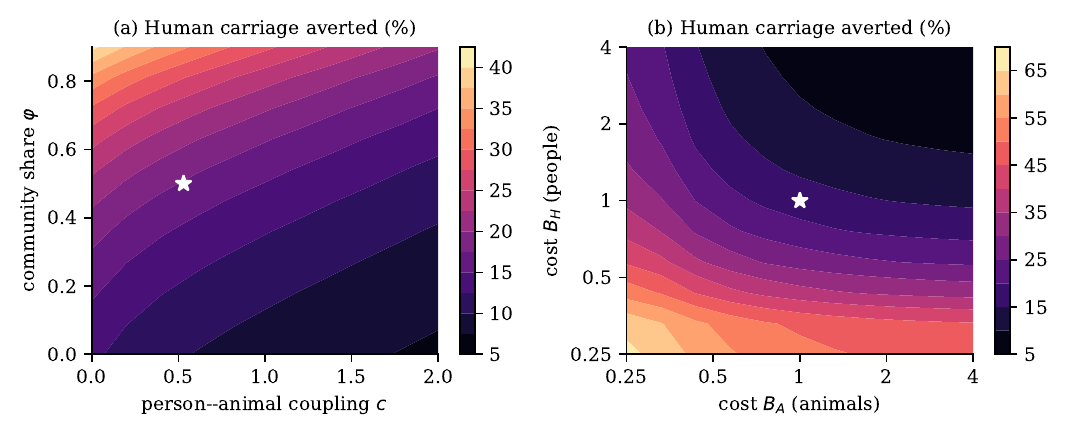}
\caption{Two-dimensional studies of the human carriage averted by the optimal strategy (\%). (a) As a function of the coupling $c$ and of the community share $\varphi$. (b) As a function of the costs of the interventions on people and on animals (logarithmic scales). Stars: baseline.}
\label{fig:twod}
\end{figure}

Figure~\ref{fig:twod} shows two pairs of parameters. In panel (a), the human carriage averted increases with the community share $\varphi$ and decreases with the coupling $c$ between person and animal, and the two interact: without coupling, raising $\varphi$ from $0$ to $0.9$ raises the proportion averted from $13\%$ to $41\%$, while with $c=2$ it raises it only from $7\%$ to $20\%$. When transmission inside the household is strong, a larger part of human carriage comes from the household's own animals, which the community interventions reach only indirectly. In panel (b), the proportion averted depends much more on the cost of the human intervention than on the cost of the animal intervention: with $B_A=1$, it falls from $52\%$ to $8\%$ as $B_H$ goes from $0.25$ to $4$, and with $B_H=1$ from $34\%$ to $14\%$ as $B_A$ goes over the same range. Moreover, each optimal intervention depends almost only on its own cost: the mean of $u_H^*$ is $0.82$, $0.22$ and $0.05$ for $B_H=0.25$, $1$ and $4$, whatever $B_A$, and similarly for $u_A^*$. In practice, each service can set the level of its own intervention from its own cost.

\subsection{Optimal controls in three dimensions}

\begin{figure}[t]
\centering
\includegraphics[width=\textwidth]{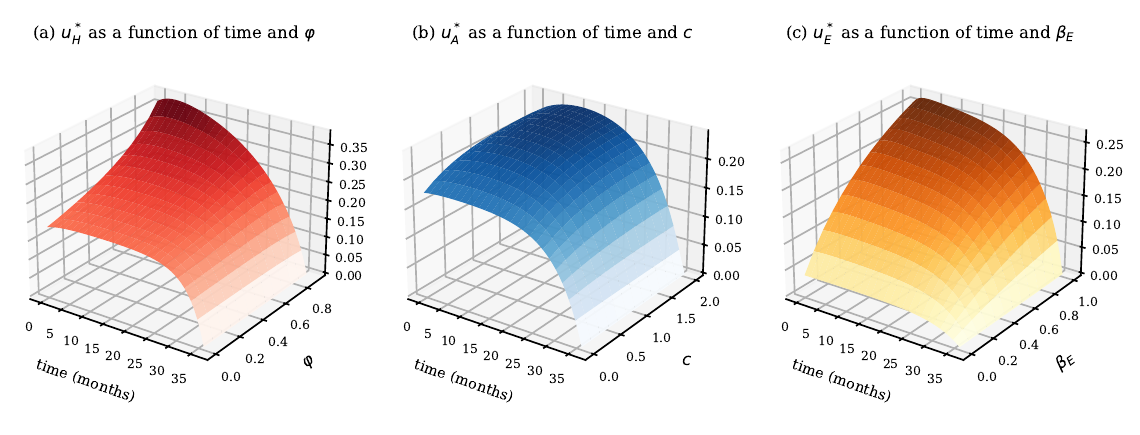}
\caption{Optimal interventions as functions of time and of one parameter: (a) $u_H^*$ and the community share $\varphi$; (b) $u_A^*$ and the coupling $c$; (c) $u_E^*$ and the rate $\beta_E$ from the soil to people and animals. Vertical axes: optimal intervention.}
\label{fig:threed}
\end{figure}

Figure~\ref{fig:threed} shows each optimal intervention as a function of time and of the parameter that affects it most. The intervention on people grows with the community share $\varphi$, at all times, since the community route is the one it acts on. The intervention on animals grows slightly with the coupling $c$, because animals then feed human carriage more directly. The intervention on the soil is small when the soil does not transmit ($\beta_E=0$; mean $0.03$, a quarter of its baseline value, since it still reduces the acquisition of the soil) and grows steadily with $\beta_E$. In all cases the shape in time is the same: a level close to the maximum for about a year and a half, then a decrease to zero at the horizon.

\begin{figure}[t]
\centering
\includegraphics[width=\textwidth]{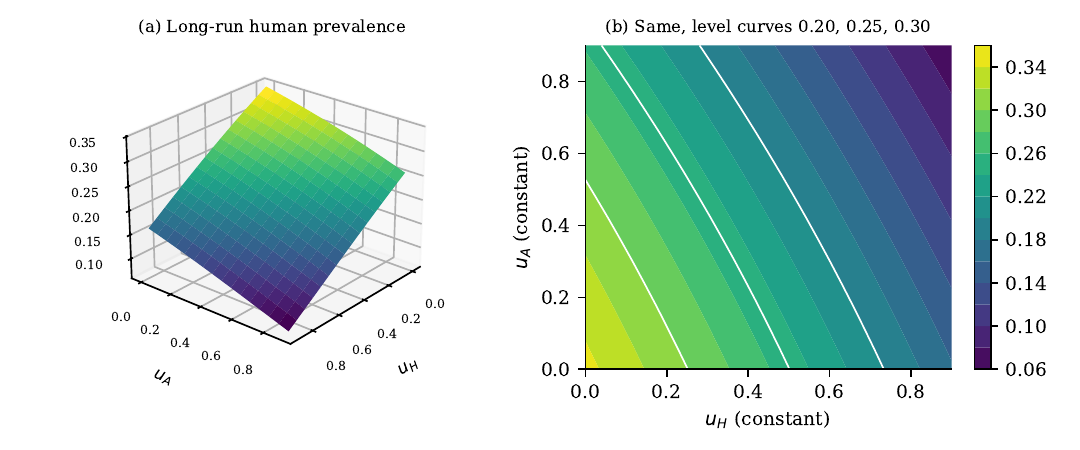}
\caption{Long-run human prevalence under constant interventions on people and animals ($u_E=0$): (a) surface and (b) level curves.}
\label{fig:equilibrium}
\end{figure}

Figure~\ref{fig:equilibrium} illustrates Proposition~\ref{prop:equilibrium}: the long-run human prevalence under constant interventions on people and animals, starting from the baseline. It decreases in both interventions, faster in $u_H$, from $0.35$ without intervention to $0.16$ with $u_H=0.9$ alone, $0.26$ with $u_A=0.9$ alone, and $0.06$ with both. For each pair of interventions on a grid of $16\times16$ values, solutions started from five different initial laws (the baseline, no node carrying the trait, all nodes carrying it, and two random laws) reached the same equilibrium to within $10^{-9}$, in agreement with Theorem~\ref{thm:scalar} (Section~\ref{sec:eqnum}).

\subsection{Equilibria and the long-run problem}\label{sec:eqnum}

\begin{figure}[t]
\centering
\includegraphics[width=\textwidth]{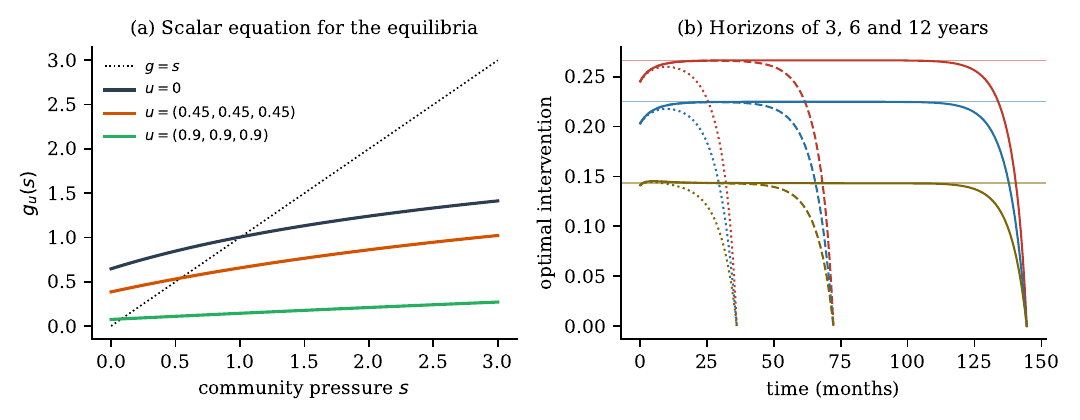}
\caption{(a) The function $g_u$ of Theorem~\ref{thm:scalar} for three constant interventions; the equilibria are its crossings with the diagonal (dotted). (b) Optimal interventions for horizons of $3$ (dotted), $6$ (dashed) and $12$ years (solid); thin horizontal lines: best constant intervention for the long-run cost. Red: people; blue: animals; brown: soil.}
\label{fig:scalar}
\end{figure}

Theorems~\ref{thm:scalar} and~\ref{thm:unique-eq} reduce the equilibrium to the unique root of $g_u(s)=s$. Figure~\ref{fig:scalar}a shows $g_u$ for three constant interventions: each curve crosses the diagonal once, and stronger interventions move the crossing, that is the equilibrium community pressure, to the left. The computations agree with the theorem: on a grid of $7\times7\times3$ constant interventions with the baseline parameters, for the $250$ parameter sets of the global sensitivity study, and for $300$ further random sets with very strong feedback ($\varphi$ between $0.9$ and $0.999$, transmission rates up to $6$), the equation always had a single root, although the slope of $g_u$ reached $6.6$ away from the root. In $400$ random sets of rates, outside rates and interventions, the inequality of Lemma~\ref{lem:subhom} held with a smallest margin of $0.0008$, and $g_u(s)/s$ was decreasing (\texttt{check\_subhom.py}).

The long-run optimum, that is the best constant intervention for the long-run cost $\bar J$, is $(u_H,u_A,u_E)=(0.267,\,0.225,\,0.143)$. It lowers the long-run prevalences from $(0.35,0.36,0.30)$ to $(0.27,0.28,0.25)$ and the long-run cost per unit of time by $10\%$. Figure~\ref{fig:scalar}b shows the optimal controls for horizons of $3$, $6$ and $12$ years: as the horizon grows, they stay longer on a plateau, and the plateau approaches the static optimum; with a horizon of $12$ years, the controls at mid-horizon are $(0.2665,\,0.2250,\,0.1434)$, equal to the static optimum to three digits. The optimal control thus has the turnpike shape known for many optimal control problems \citep{porretta2013,trelat2015}: a short adjustment at the start, a long plateau at the long-run optimum, and a decrease at the end that only reflects the finite horizon. We distinguish three levels: the analytical results (Theorems~\ref{thm:turnpike}, \ref{thm:lqturnpike} and~\ref{thm:dissip}), the numerical verification of their hypotheses for our parameters, and the numerical observation of the turnpike for the nonlinear problem.

\emph{Verification of the hypotheses of Theorem~\ref{thm:lqturnpike}.} The assumptions of Theorem~\ref{thm:lqturnpike} hold at the static optimum, where $p_s(x)\ge0.05$ for all $x$ and $u_s$ is interior: $H_{vv}$ is the identity, since $B_k=1$; the matrix $\mathbb H$ has no eigenvalue on the imaginary axis, the real parts closest to it being $\pm0.574$; and the stable and unstable subspaces are graphs over $\xi$ and $\mu$ (smallest singular values of the normalised blocks $X_s$ and $M_u$: $0.53$ and $0.61$). Moreover, $A$ is invertible and stable, with eigenvalues of real parts between $-5.28$ and $-0.61$, the pair $(A,B_u)$ is controllable (the Kalman matrix has rank $7$), and the reduced Hessian of the static problem on the set of equilibria has eigenvalues $0.76, 0.94, 1.03$, so the static optimum is a strict local minimum. 

\emph{Observation for the nonlinear problem.} For the nonlinear problem with a horizon of $12$ years, the distance $d$ of the optimal pair to the static optimum decreased exponentially at the rate $0.574$ per time unit at the start and increased at the rate $0.574$ at the end, in agreement with the smallest absolute real part $0.574$ of the eigenvalues of $\mathbb H$, the rate allowed by Theorem~\ref{thm:lqturnpike}; at mid-horizon it was $10^{-5}$. In calendar time, the distance is divided by $e$ every $1/\sigma\approx1.7$ time units (about six months), and by ten in about $4$ time units (about $14$ months); the final layer has the same width. 

\emph{Verification of Assumption~\ref{ass:dissip}.} With $\varphi=0$ and the other baseline rates, the constants of Theorem~\ref{thm:dissip} depend on the costs through the static optimum; with $B_k=b$ for all $k$, $4C_2C_4$ equals $1.89$ for $b=1$ and $1.74$ for $b=1.9$, so the theorem guarantees the turnpike property for $b\ge1.9$ (the static optimum stays interior). With $\varphi=0.1$ and $b=1$ the requirement is $B_k>6.5$, and for $\varphi\ge0.3$ our bound on $c_K$ is no longer positive. These conditions are sufficient, not necessary. For the baseline parameters we therefore checked Assumption~\ref{ass:dissip} directly with the storage $\Sigma_\varepsilon$, $\varepsilon=4$. In the coordinates $(\xi,\delta u)$, the function $\ell(p,u)-\ell_s+\nabla\Sigma_\varepsilon(p)\cdot f(p,u)-\kappa(|\xi|^2+|\delta u|^2)$ is a polynomial of degree four, and $K\times[0,\bar u]^3$ is a polytope described by the inequalities $p\ge0$, $p_{\min}(\bar u\mathbf 1)(\Gamma)\le p(\Gamma)\le\bar p(\Gamma)$ for the $18$ nonempty proper increasing sets $\Gamma$, and $0\le u_k\le\bar u$. We wrote this polynomial as a sum of squares plus a combination of the constraints with sum-of-squares or nonnegative multipliers (products of pairs of linear constraints, and the quadratic constraints $u_k(\bar u-u_k)\ge0$, $p(x)(1-p(x))\ge0$ and $|\xi|^2+|\delta u|^2\le2+3\bar u^2$), and maximised $\kappa$ by semidefinite programming (\texttt{sos\_dissip.py}, solver Clarabel). The largest certified value is $\kappa=0.287$, which equals the smallest ratio found by local minimisation from $80$ starting points; the largest difference between the coefficients of the two sides of the polynomial identity is $3\times10^{-8}$ and the smallest eigenvalues of the Gram matrices are above $-2\times10^{-9}$. Certificates of the same kind exist with $\varphi=0$ ($\varepsilon=2$, $\kappa=0.31$), $\varphi=0.25$ ($\varepsilon=4$, $\kappa=0.32$) and $\varphi=0.75$ ($\varepsilon=8$, $\kappa=0.13$), with coefficient differences below $6\times10^{-8}$ and Gram eigenvalues above $-3\times10^{-8}$, that is over the whole plausible range of the community share; with $\varphi=0.9$ we did not find one with storages of this form. This is a computer-assisted numerical certification, not an exact proof: the certificate is computed in floating point by an interior-point solver, and an exact certificate would require rounding the Gram matrices to rational numbers and checking the identity in exact arithmetic, which the margin $\kappa=0.287$ should allow but which we have not done. Up to the accuracy of the solver, Assumption~\ref{ass:dissip} therefore holds for our parameters, and Theorem~\ref{thm:turnpike} applies. With $C_0=0.20$ and the oscillation of $\Sigma_4$ on $K$ equal to $1.23$, it gives $|V_T-T\ell_s|\le1.23$ and $J_T(u_s)-V_T\le1.43$ for every horizon: over $12$ years the relative loss of the constant long-run optimum, $(J_T(u_s)-V_T)/J_T(u_s)$, is below $7\%$, and over longer horizons it tends to zero. The bound on the time spent away from the static optimum, $C_1/(\kappa r^2)$, is far from sharp (about $20$ time units at distance $r=0.5$), whereas the computed optimal control leaves the neighbourhood only during the two layers of Theorem~\ref{thm:lqturnpike}. For a programme of six years or more, the long-run optimum, which is computed from Theorem~\ref{thm:scalar} by a scalar root and a minimisation over three numbers, is a simple policy that is close to optimal. It should not be confused with the best constant intervention over a given horizon, $(0.22,0.18,0.12)$ for three years, which is smaller because the horizon is short.

\subsection{Comparison with a compartmental model}

To see what the household structure changes, we compared our model with the compartmental model that a population-level study would use: three populations, people, animals and soil, with prevalences $m_k$ that follow
\[
\dot m_k=(1-m_k)\Big[(1-u_k)a_k^{c}(m)+\sum_{j\neq k}\theta_j(u)\beta_{jk}m_j\Big]-\gamma_km_k ,
\]
with the same rates $\beta$, $\gamma$, $\varphi$ and $\omega$, and outside rates $\alpha^c$ chosen so that the baseline prevalences are the same, $\bar m$. This model needs smaller outside rates, $\alpha^c=(0.25,\,0.28,\,0.09)$ instead of $(0.31,\,0.34,\,0.15)$: it lets every person meet the average animal, while in the household model the carriers are grouped in the same households, so that a person who does not carry the trait lives less often than average with an animal that does. The compartmental model thus attributes more carriage to transmission and less to the outside routes on which the interventions act.

The consequences are shown in Figure~\ref{fig:zones}a. We computed the optimal intervention of the compartmental model (by direct minimisation over controls constant on $40$ intervals). It is weaker on people and animals and stronger on the soil (mean $(0.17,\,0.14,\,0.13)$ against $(0.22,\,0.18,\,0.12)$), and the compartmental model predicts that it averts $11.7\%$ of human carriage. Applied to the household model, it averts $13.8\%$ (so the compartmental model underestimates the effect of its own policy by $15\%$), against $17.3\%$ for the household optimum, and it achieves $95\%$ of the reduction of $J$ that is possible. For the household optimum, the compartmental model predicts $14.4\%$ instead of $17.3\%$, an underestimate of $17\%$. With the same transmission rates, ignoring the household therefore underestimates the benefit of the interventions by about a sixth, and leads to a policy that puts less effort on people and animals and loses about $5\%$ of the achievable benefit. A compartmental model fitted to population data would have its own transmission rates, and the size of these differences would change; we expect the direction of the bias, which comes from the clustering of carriers in households, to remain, but we have shown it only for the calibration used here.

\subsection{Two zones of the town}

The households of the survey came from the formal (\emph{loti}, $37$ households) and informal (\emph{non loti}, $22$ households) zones of Tenkodogo. Pooling the five resistance traits, the observed prevalences were $0.375$, $0.344$ and $0.525$ in people, animals and soil in the formal zone, and $0.250$, $0.347$ and $0.333$ in the informal zone. We built a model with two types of households in these proportions, sharing one community, with the within-household rates of Table~\ref{tab:params}, and with outside rates chosen so that the stationary prevalences of each type are proportional to these observed values and average to $\bar m$. This gives $\alpha^{\rm loti}=(0.42,\,0.29,\,0.23)$ and $\alpha^{\text{non loti}}=(0.16,\,0.41,\,0.04)$: in the formal zone people and soil acquire more from outside, in the informal zone the animals do. These zone-specific values rest on few isolates and are illustrative.

The targeted strategy (Figure~\ref{fig:zones}b) puts the intervention on people and soil mainly in the formal zone ($u_H=0.27$ and $u_E=0.15$ on average, against $0.13$ and $0.06$ in the informal zone), and the intervention on animals mainly in the informal zone ($0.25$ against $0.15$). It averts $18.4\%$ of human carriage and reduces $J$ by $8.1\%$, against $16.5\%$ and $7.4\%$ for the best uniform strategy, for a cost of the interventions about $12\%$ higher. The gain is not shared equally: the targeted strategy averts $19.3\%$ of human carriage in the formal zone and $16.0\%$ in the informal zone, against $16.1\%$ and $17.3\%$ for the uniform strategy. Targeting by zone improves the overall result when the exposure of the zones differs, and the model says which intervention to put where, but it lowers the protection of people in the zone where human exposure is smaller; the choice between the two strategies involves equity as well as efficiency.

\subsection{Robustness to the community share}

\begin{figure}[t]
\centering
\includegraphics[width=\textwidth]{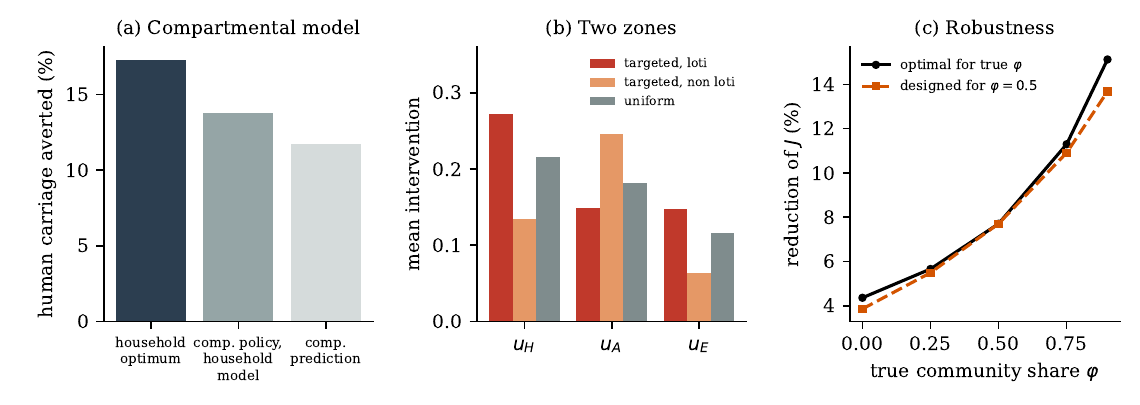}
\caption{(a) Human carriage averted: household-model optimum; optimum of the compartmental model applied to the household model; and the benefit predicted by the compartmental model for its own optimum. (b) Mean optimal interventions in the formal (\emph{loti}) and informal (\emph{non loti}) zones, targeted strategy, and best uniform strategy. (c) Reduction of the total cost $J$ with the strategy optimal for the true community share $\varphi$, and with the strategy designed for $\varphi=0.5$.}
\label{fig:zones}
\end{figure}

The community share $\varphi$ is the parameter that matters most and the least known. We applied the optimal strategy computed with $\varphi=0.5$ to models with other values of $\varphi$, and compared it with the strategy that would be optimal if $\varphi$ were known (Figure~\ref{fig:zones}c). The strategy designed with $\varphi=0.5$ achieves $88\%$ of the best possible reduction of $J$ when $\varphi=0$, $97\%$ when $\varphi=0.25$ or $0.75$, and $90\%$ when $\varphi=0.9$. A wrong value of $\varphi$ changes the predicted benefit a lot, but the strategy itself remains close to optimal.

\subsection{Robust strategy and value of information}

We took nine scenarios, combining $\varphi\in\{0.25,\,0.5,\,0.75\}$, the plausible range of the community share, with $c\in\{0.2,\,0.53,\,1.0\}$, which covers most of the confidence interval of the coupling, with equal weights, and computed the strategy that minimises the expected cost. It is very close to the baseline strategy: its mean interventions are $(0.23, 0.18, 0.12)$ against $(0.22, 0.18, 0.12)$. In every scenario it achieves at least $94\%$ of the reduction of $J$ that the strategy optimal for that scenario would achieve, and the baseline strategy at least $95\%$ (Table~\ref{tab:robust}). The expected value of perfect information is $2.3\%$ of the expected achievable reduction of $J$.

This has a practical consequence. The uncertain parameters change the predicted benefit a lot, from $11\%$ to $28\%$ of human carriage averted in these scenarios, but they change the best strategy very little. A new survey to estimate $\varphi$ and $c$ would improve the prediction of the benefit, but it would improve the decision by only about $2\%$. A health service can therefore choose its interventions now, with the present data, and use later data to evaluate them.

\begin{table}[t]
\centering
\caption{Robust strategy. For each scenario: human carriage averted by the strategy optimal for that scenario and by the robust strategy, and share of the best achievable reduction of $J$ obtained by the robust strategy and by the baseline strategy (computed with $\varphi=0.5$, $c=0.53$).}
\label{tab:robust}
\small
\begin{tabular}{cccccc}
\toprule
$\varphi$ & $c$ & \multicolumn{2}{c}{Human carriage averted (\%)} & \multicolumn{2}{c}{Share of best reduction of $J$ (\%)}\\
 & & scenario optimum & robust & robust & baseline strategy\\
\midrule
0.25 & 0.20 & 14.7 & 17.4 & 94.3 & 95.4\\
0.25 & 0.53 & 12.9 & 15.3 & 96.1 & 97.0\\
0.25 & 1.00 & 11.2 & 13.3 & 96.9 & 97.7\\
0.50 & 0.20 & 19.7 & 20.1 & 99.6 & 99.7\\
0.50 & 0.53 & 17.3 & 17.6 & 100.0 & 100.0\\
0.50 & 1.00 & 14.9 & 15.1 & 99.8 & 99.8\\
0.75 & 0.20 & 28.2 & 24.3 & 97.8 & 97.3\\
0.75 & 0.53 & 24.8 & 20.9 & 97.2 & 96.6\\
0.75 & 1.00 & 21.3 & 17.5 & 96.7 & 96.0\\
\bottomrule
\end{tabular}
\end{table}

%======================================================================
\section{Practical guide}\label{sec:guide}
%======================================================================

Table~\ref{tab:guide} summarises the results for public health services.

\begin{table}[t]
\centering
\caption{Practical guide.}
\label{tab:guide}
\small
\begin{tabular}{p{4.3cm}p{10.2cm}}
\toprule
Question & Answer from the model\\
\midrule
Can an intervention make things worse? & No. In this model every intervention, on people, animals or soil, lowers carriage everywhere and at all times (Corollary~\ref{cor:norebound}). \\[3pt]
Should one intervene on one route only? & No, under a quadratic cost of effort: all three interventions are part of the optimal strategy at every time (Theorem~\ref{thm:positive}). With the Tenkodogo parameters they also avert together more than the sum of their separate effects, because of the transmission inside households (Table~\ref{tab:strategies}). \\[3pt]
Where should most effort go? & With the parameters of Tenkodogo and equal costs, to people (prescription-only dispensing), then animals (veterinary access), then soil (treatment of sludge and manure). This order is not general: each intervention grows with the acquisition rate of its route and decreases with its own cost, so the order must be recomputed with local costs (Section~\ref{sec:numerics}). \\[3pt]
How should effort vary in time? & For a three-year programme, close to constant for the first year and a half, then decreasing. For programmes of six years or more, close to the long-run optimum, computed from a scalar equation (Theorem~\ref{thm:scalar}), except in the last year. \\[3pt]
Should the interventions differ between zones? & It can help when exposure differs: in Tenkodogo the model puts the human and soil interventions mainly in the formal zone and the animal intervention in the informal zone. This averts more carriage overall than a uniform strategy, but less in the informal zone, so equity must be weighed. \\[3pt]
Can a simpler model be used? & With the same transmission rates, a model of three populations without households underestimates the benefit by about a sixth and recommends less effort on people and animals. \\[3pt]
Is a new survey needed before acting? & Not for choosing the strategy: the robust strategy is within a few percent of the best one in all plausible scenarios, and the value of perfect information on $\varphi$ and $c$ is about $2\%$ (Table~\ref{tab:robust}). New data mainly improve the prediction of the benefit. \\[3pt]
What needs to be measured first? & To predict the benefit: the share of acquisition due to the community and the costs of the interventions, which drive the result; then the duration of carriage in animals. The soil parameters matter little for human carriage, but they determine the effort on the soil. \\[3pt]
Does the model apply to a village? & Only roughly: below about a thousand households, chance fluctuations are not small compared with the effect of the interventions (Figure~\ref{fig:mf}). \\
\bottomrule
\end{tabular}
\end{table}

%======================================================================
\section{Discussion}\label{sec:discussion}
%======================================================================

The paper combines a household model, which keeps the joint state of people, animals and soil, with community feedback, and studies its optimal control with complete proofs. The household structure matters for the results: the synergy between interventions comes from the transmission inside the household, the monotonicity results rest on the attractiveness of the household chain, and a compartmental model with the same transmission rates underestimates the benefit of the interventions because it ignores the clustering of carriers in households. For long programmes, the scalar equation of Theorem~\ref{thm:scalar} gives a long-run optimum that is close to the optimal strategy, which makes the results easy to use; the turnpike theorems explain why, and the linearised problem predicts that the optimal strategy approaches this optimum with an e-folding time of about six months. Strict dissipativity, which gives the turnpike property, holds under the explicit condition of Theorem~\ref{thm:dissip}, and we have verified it by a computer-assisted numerical certification for our parameters over the plausible range of the community share; a proof for all parameters, or the construction of a storage that works for very strong community feedback, remains open. The community feedback matters as much: at the assumed community share $\varphi=0.5$, a model without it underestimates the benefit of the optimal strategy by about two fifths; the size of this effect depends on $\varphi$, which is not identified by the data.

The mathematical results extend to more nodes. Theorems~\ref{thm:mf},~\ref{thm:comparison} and~\ref{thm:positive} hold, with constants that depend on the number $2^n$ of household states, for any attractive single-flip dynamics with $n$ nodes in which the community term is nondecreasing and the interventions only reduce acquisition rates; this covers several people and animals per household, or water and food as further nodes. The comparison theorem is the main tool: the mean-field system is monotone for the stochastic order of household states, that is for the cone of the stochastic order, rather than for the componentwise order of $p$ \citep{smith1995}.

The study has limits. The interventions are fixed in advance for each type of household; strategies that react to the observed state of individual households would require the control of the finite system or of a mean-field game. The costs $B_k$ are relative and not in money, and the weights $w$ are a policy choice; the sensitivity study shows that the qualitative results (all interventions used, synergy, the main drivers) hold over wide ranges, but the levels of effort, and even the order of the three interventions, depend on these values. The links with the soil and the community share $\varphi$ cannot be estimated from a survey at one time. The community share drives the predicted benefit, but the robust analysis shows that it changes the best strategy little, so that the strategy can be chosen with the present data; estimating $\varphi$, for example by sampling households at two times in places with different interventions, would mainly sharpen the prediction. In the comparison with a compartmental model we kept the transmission rates of the household model, because a population model fitted only to the prevalences of the three nodes cannot identify its own transmission rates. The model has two states per node, so it does not distinguish strains or resistance acquired under treatment, and it assumes that the durations of carriage are not changed by the interventions. Finally, the two-zone model takes into account only one source of difference between households; differences within a zone would add a household-level random effect to the acquisition rates, which can be approximated by further types of households, with the same results.

%======================================================================
\section{Conclusion and remarks}\label{sec:conclusion}
%======================================================================

A household model of antimicrobial resistance linked through the community has a mean-field limit that can be analysed and controlled: it is monotone, its equilibria are given by a scalar equation, and its optimal control exists and can be characterised. Interventions on people, animals and soil never increase carriage and, with a quadratic cost, are all part of the optimal strategy; with the parameters of Tenkodogo they act together more than separately. In the baseline scenario, a combined programme of three years averts about a sixth of human carriage under moderate costs (between $11\%$ and $28\%$ over the plausible scenarios), on long horizons the best constant intervention is nearly optimal, and the choice of the strategy is robust to the parameters that a single survey cannot estimate. The household structure matters for these conclusions, and the share of acquisition due to the community is the quantity that most needs to be measured. The same framework can include more nodes, more types of households and costs in money when they are known.

\section*{Acknowledgements}
We thank the households of Tenkodogo and the health, veterinary and environment services of the Centre-Est region.

\section*{Code availability}
The Python code that computes all results and figures of the paper is provided as supplementary material and will be deposited in a public archive on acceptance.

\section*{Conflict of interest}
None declared.

\bibliographystyle{plainnat}
\bibliography{refs}

\begin{thebibliography}{35}
\providecommand{\natexlab}[1]{#1}
\providecommand{\url}[1]{\texttt{#1}}
\expandafter\ifx\csname urlstyle\endcsname\relax
  \providecommand{\doi}[1]{doi: #1}\else
  \providecommand{\doi}{doi: \begingroup \urlstyle{rm}\Url}\fi

\bibitem[Aizouk et~al.(2025)Aizouk, Sibidou, Valia, Ingelbeen, Campbell,
  Kouanda, Welgo, Diagne, Riems, Hardy, Wouters, Meudec, van~der Sande,
  Halidou, Cooper, and van Kleef]{aizouk2025}
R.~Aizouk, Y.~Sibidou, D.~Valia, B.~Ingelbeen, L.~Campbell, J.~S. Kouanda,
  A.~Welgo, P.~M. Diagne, B.~Riems, L.~Hardy, E.~Wouters, M.~Meudec, M.~A.~B.
  van~der Sande, T.~Halidou, B.~S. Cooper, and E.~van Kleef.
\newblock Targeting community-level drivers of antimicrobial resistance in
  sub-{S}aharan {A}frica: the effect of a community-based intervention bundle
  on household transmission of extended spectrum beta-lactamase-producing
  \emph{E. coli} in rural {B}urkina {F}aso -- a cluster randomised trial.
\newblock medRxiv preprint, doi:10.64898/2025.12.15.25342269, 2025.

\bibitem[{Antimicrobial Resistance Collaborators}(2022)]{amr2022}
{Antimicrobial Resistance Collaborators}.
\newblock Global burden of bacterial antimicrobial resistance in 2019: a
  systematic analysis.
\newblock \emph{The Lancet}, 399\penalty0 (10325):\penalty0 629--655, 2022.
\newblock \doi{10.1016/S0140-6736(21)02724-0}.

\bibitem[Ball et~al.(1997)Ball, Mollison, and Scalia-Tomba]{ball1997}
Frank Ball, Denis Mollison, and Gianpaolo Scalia-Tomba.
\newblock Epidemics with two levels of mixing.
\newblock \emph{The Annals of Applied Probability}, 7\penalty0 (1):\penalty0
  46--89, 1997.
\newblock \doi{10.1214/aoap/1034625252}.

\bibitem[Bauer et~al.(1966)Bauer, Kirby, Sherris, and Turck]{bauer1966}
A.~W. Bauer, W.~M.~M. Kirby, J.~C. Sherris, and M.~Turck.
\newblock Antibiotic susceptibility testing by a standardized single disk
  method.
\newblock \emph{American Journal of Clinical Pathology}, 45\penalty0
  (4):\penalty0 493--496, 1966.
\newblock \doi{10.1093/ajcp/45.4_ts.493}.

\bibitem[Bonhoeffer et~al.(1997)Bonhoeffer, Lipsitch, and
  Levin]{bonhoeffer1997}
Sebastian Bonhoeffer, Marc Lipsitch, and Bruce~R. Levin.
\newblock Evaluating treatment protocols to prevent antibiotic resistance.
\newblock \emph{Proceedings of the National Academy of Sciences of the USA},
  94\penalty0 (22):\penalty0 12106--12111, 1997.
\newblock \doi{10.1073/pnas.94.22.12106}.

\bibitem[Booton et~al.(2021)Booton, Meeyai, Alhusein, Buller, Feil, Lambert,
  Mongkolsuk, Pitchforth, Reyher, Sakcamduang, Satayavivad, Singer,
  Sringernyuang, Thamlikitkul, Vass, {OH-DART Study Group}, Avison, and
  Turner]{booton2021}
Ross~D. Booton, Aronrag Meeyai, Nour Alhusein, Henry Buller, Edward Feil, Helen
  Lambert, Skorn Mongkolsuk, Emma Pitchforth, Kristen~K. Reyher, Walasinee
  Sakcamduang, Jutamaad Satayavivad, Andrew~C. Singer, Luechai Sringernyuang,
  Visanu Thamlikitkul, Lucy Vass, {OH-DART Study Group}, Matthew~B. Avison, and
  Katherine M.~E. Turner.
\newblock One {H}ealth drivers of antibacterial resistance: quantifying the
  relative impacts of human, animal and environmental use and transmission.
\newblock \emph{One Health}, 12:\penalty0 100220, 2021.
\newblock \doi{10.1016/j.onehlt.2021.100220}.

\bibitem[Chilanga et~al.(2025)Chilanga, Kasozi, Mazeri, Paterson, and
  Muwonge]{chilanga2025}
Frank Chilanga, Keneth~I. Kasozi, Stella Mazeri, Gavin~K. Paterson, and Adrian
  Muwonge.
\newblock A systematic review of antimicrobial resistance transmission
  inferences at the human-livestock interface in {A}frica.
\newblock \emph{npj Antimicrobials and Resistance}, 3:\penalty0 58, 2025.
\newblock \doi{10.1038/s44259-025-00126-y}.

\bibitem[Ethier and Kurtz(1986)]{ethier1986}
Stewart~N. Ethier and Thomas~G. Kurtz.
\newblock \emph{Markov Processes: Characterization and Convergence}.
\newblock Wiley, New York, 1986.

\bibitem[Faulwasser et~al.(2017)Faulwasser, Korda, Jones, and
  Bonvin]{faulwasser2017}
Timm Faulwasser, Milan Korda, Colin~N. Jones, and Dominique Bonvin.
\newblock On turnpike and dissipativity properties of continuous-time optimal
  control problems.
\newblock \emph{Automatica}, 81:\penalty0 297--304, 2017.
\newblock \doi{10.1016/j.automatica.2017.03.012}.

\bibitem[Fister et~al.(1998)Fister, Lenhart, and McNally]{fister1998}
K.~Renee Fister, Suzanne Lenhart, and Joseph~Scott McNally.
\newblock Optimizing chemotherapy in an {HIV} model.
\newblock \emph{Electronic Journal of Differential Equations}, 1998\penalty0
  (32):\penalty0 1--12, 1998.

\bibitem[Fleming and Rishel(1975)]{fleming1975}
Wendell~H. Fleming and Raymond~W. Rishel.
\newblock \emph{Deterministic and Stochastic Optimal Control}.
\newblock Springer, New York, 1975.

\bibitem[Gr{\"u}ne and M{\"u}ller(2016)]{grune2016}
Lars Gr{\"u}ne and Matthias~A. M{\"u}ller.
\newblock On the relation between strict dissipativity and turnpike properties.
\newblock \emph{Systems \& Control Letters}, 90:\penalty0 45--53, 2016.

\bibitem[Haverkate et~al.(2017)Haverkate, Platteel, Fluit, Cohen~Stuart,
  Leverstein-van Hall, Thijsen, Scharringa, Kloosterman, Bonten, and
  Bootsma]{haverkate2017}
Manon~R. Haverkate, Tamara~N. Platteel, A.~C. Fluit, James~W. Cohen~Stuart,
  Maurine~A. Leverstein-van Hall, Steven F.~T. Thijsen, J.~Scharringa, Fieke
  R.~C. Kloosterman, Marc J.~M. Bonten, and Martin C.~J. Bootsma.
\newblock Quantifying within-household transmission of extended-spectrum
  $\beta$-lactamase-producing bacteria.
\newblock \emph{Clinical Microbiology and Infection}, 23\penalty0 (1):\penalty0
  46.e1--46.e7, 2017.
\newblock \doi{10.1016/j.cmi.2016.08.021}.

\bibitem[House and Keeling(2008)]{house2008}
Thomas House and Matt~J. Keeling.
\newblock Deterministic epidemic models with explicit household structure.
\newblock \emph{Mathematical Biosciences}, 213\penalty0 (1):\penalty0 29--39,
  2008.
\newblock \doi{10.1016/j.mbs.2008.01.011}.

\bibitem[Kamae et~al.(1977)Kamae, Krengel, and O'Brien]{kamae1977}
Teturo Kamae, Ulrich Krengel, and George~L. O'Brien.
\newblock Stochastic inequalities on partially ordered spaces.
\newblock \emph{The Annals of Probability}, 5\penalty0 (6):\penalty0 899--912,
  1977.
\newblock \doi{10.1214/aop/1176995659}.

\bibitem[Kurtz(1970)]{kurtz1970}
Thomas~G. Kurtz.
\newblock Solutions of ordinary differential equations as limits of pure jump
  {M}arkov processes.
\newblock \emph{Journal of Applied Probability}, 7\penalty0 (1):\penalty0
  49--58, 1970.
\newblock \doi{10.2307/3212147}.

\bibitem[Lasserre(2001)]{lasserre2001}
Jean~B. Lasserre.
\newblock Global optimization with polynomials and the problem of moments.
\newblock \emph{SIAM Journal on Optimization}, 11\penalty0 (3):\penalty0
  796--817, 2001.
\newblock \doi{10.1137/S1052623400366802}.

\bibitem[Lenhart and Workman(2007)]{lenhart2007}
Suzanne Lenhart and John~T. Workman.
\newblock \emph{Optimal Control Applied to Biological Models}.
\newblock Chapman and Hall/CRC, Boca Raton, 2007.

\bibitem[Liggett(1985)]{liggett1985}
Thomas~M. Liggett.
\newblock \emph{Interacting Particle Systems}, volume 276 of \emph{Grundlehren
  der mathematischen Wissenschaften}.
\newblock Springer, New York, 1985.
\newblock \doi{10.1007/b138374}.
\newblock Reprinted in Classics in Mathematics, Springer, Berlin, 2005.

\bibitem[Lindsay(1988)]{lindsay1988}
Bruce~G. Lindsay.
\newblock Composite likelihood methods.
\newblock In N.~U. Prabhu, editor, \emph{Statistical Inference from Stochastic
  Processes}, volume~80 of \emph{Contemporary Mathematics}. American
  Mathematical Society, Providence, RI, 1988.

\bibitem[Lipsitch et~al.(2000)Lipsitch, Bergstrom, and Levin]{lipsitch2000}
Marc Lipsitch, Carl~T. Bergstrom, and Bruce~R. Levin.
\newblock The epidemiology of antibiotic resistance in hospitals: paradoxes and
  prescriptions.
\newblock \emph{Proceedings of the National Academy of Sciences of the USA},
  97\penalty0 (4):\penalty0 1938--1943, 2000.
\newblock \doi{10.1073/pnas.97.4.1938}.

\bibitem[Marino et~al.(2008)Marino, Hogue, Ray, and Kirschner]{marino2008}
Simeone Marino, Ian~B. Hogue, Christian~J. Ray, and Denise~E. Kirschner.
\newblock A methodology for performing global uncertainty and sensitivity
  analysis in systems biology.
\newblock \emph{Journal of Theoretical Biology}, 254\penalty0 (1):\penalty0
  178--196, 2008.
\newblock \doi{10.1016/j.jtbi.2008.04.011}.

\bibitem[Massey(1987)]{massey1987}
William~A. Massey.
\newblock Stochastic orderings for {M}arkov processes on partially ordered
  spaces.
\newblock \emph{Mathematics of Operations Research}, 12\penalty0 (2):\penalty0
  350--367, 1987.
\newblock \doi{10.1287/moor.12.2.350}.

\bibitem[McKay et~al.(1979)McKay, Beckman, and Conover]{mckay1979}
Michael~D. McKay, Richard~J. Beckman, and William~J. Conover.
\newblock A comparison of three methods for selecting values of input variables
  in the analysis of output from a computer code.
\newblock \emph{Technometrics}, 21\penalty0 (2):\penalty0 239--245, 1979.
\newblock \doi{10.2307/1268522}.

\bibitem[Mughini-Gras et~al.(2019)Mughini-Gras, Dorado-Garc{\'i}a, van
  Duijkeren, van~den Bunt, Dierikx, Bonten, Bootsma, Schmitt, Hald, Evers,
  de~Koeijer, van Pelt, Franz, Mevius, Heederik, and {ESBL Attribution
  Consortium}]{mughinigras2019}
Lapo Mughini-Gras, Alejandro Dorado-Garc{\'i}a, Engeline van Duijkeren, Gerrita
  van~den Bunt, Cindy~M. Dierikx, Marc J.~M. Bonten, Martin C.~J. Bootsma,
  Heike Schmitt, Tine Hald, Eric~G. Evers, Aline de~Koeijer, Wilfrid van Pelt,
  Eelco Franz, Dik~J. Mevius, Dick J.~J. Heederik, and {ESBL Attribution
  Consortium}.
\newblock Attributable sources of community-acquired carriage of
  \emph{Escherichia coli} containing $\beta$-lactam antibiotic resistance
  genes: a population-based modelling study.
\newblock \emph{The Lancet Planetary Health}, 3\penalty0 (8):\penalty0
  e357--e369, 2019.
\newblock \doi{10.1016/S2542-5196(19)30130-5}.

\bibitem[Parrilo(2003)]{parrilo2003}
Pablo~A. Parrilo.
\newblock Semidefinite programming relaxations for semialgebraic problems.
\newblock \emph{Mathematical Programming, Series B}, 96\penalty0 (2):\penalty0
  293--320, 2003.
\newblock \doi{10.1007/s10107-003-0387-5}.

\bibitem[Pontryagin et~al.(1962)Pontryagin, Boltyanskii, Gamkrelidze, and
  Mishchenko]{pontryagin1962}
Lev~S. Pontryagin, Vladimir~G. Boltyanskii, Revaz~V. Gamkrelidze, and
  Evgenii~F. Mishchenko.
\newblock \emph{The Mathematical Theory of Optimal Processes}.
\newblock Interscience, New York, 1962.

\bibitem[Porretta and Zuazua(2013)]{porretta2013}
Alessio Porretta and Enrique Zuazua.
\newblock Long time versus steady state optimal control.
\newblock \emph{SIAM Journal on Control and Optimization}, 51\penalty0
  (6):\penalty0 4242--4273, 2013.
\newblock \doi{10.1137/130907239}.

\bibitem[Smith(1995)]{smith1995}
Hal~L. Smith.
\newblock \emph{Monotone Dynamical Systems: An Introduction to the Theory of
  Competitive and Cooperative Systems}, volume~41 of \emph{Mathematical Surveys
  and Monographs}.
\newblock American Mathematical Society, Providence, 1995.

\bibitem[Strassen(1965)]{strassen1965}
Volker Strassen.
\newblock The existence of probability measures with given marginals.
\newblock \emph{The Annals of Mathematical Statistics}, 36\penalty0
  (2):\penalty0 423--439, 1965.
\newblock \doi{10.1214/aoms/1177700153}.

\bibitem[Tang et~al.(2017)Tang, Caffrey, N{\'o}brega, Cork, Ronksley, Barkema,
  Polachek, Ganshorn, Sharma, Kellner, and Ghali]{tang2017}
Karen~L. Tang, Niamh~P. Caffrey, Diego~B. N{\'o}brega, Susan~C. Cork, Paul~E.
  Ronksley, Herman~W. Barkema, Alicia~J. Polachek, Heather Ganshorn, Nishan
  Sharma, James~D. Kellner, and William~A. Ghali.
\newblock Restricting the use of antibiotics in food-producing animals and its
  associations with antibiotic resistance in food-producing animals and human
  beings: a systematic review and meta-analysis.
\newblock \emph{The Lancet Planetary Health}, 1\penalty0 (8):\penalty0
  e316--e327, 2017.
\newblock \doi{10.1016/S2542-5196(17)30141-9}.

\bibitem[Tr{\'e}lat and Zuazua(2015)]{trelat2015}
Emmanuel Tr{\'e}lat and Enrique Zuazua.
\newblock The turnpike property in finite-dimensional nonlinear optimal
  control.
\newblock \emph{Journal of Differential Equations}, 258\penalty0 (1):\penalty0
  81--114, 2015.
\newblock \doi{10.1016/j.jde.2014.09.005}.

\bibitem[van Bunnik and Woolhouse(2017)]{vanbunnik2017}
B.~A.~D. van Bunnik and M.~E.~J. Woolhouse.
\newblock Modelling the impact of curtailing antibiotic usage in food animals
  on antibiotic resistance in humans.
\newblock \emph{Royal Society Open Science}, 4\penalty0 (4):\penalty0 161067,
  2017.
\newblock \doi{10.1098/rsos.161067}.

\bibitem[Varin et~al.(2011)Varin, Reid, and Firth]{varin2011}
Cristiano Varin, Nancy Reid, and David Firth.
\newblock An overview of composite likelihood methods.
\newblock \emph{Statistica Sinica}, 21\penalty0 (1):\penalty0 5--42, 2011.

\bibitem[{World Health Organization}(2015)]{who2015}
{World Health Organization}.
\newblock Global action plan on antimicrobial resistance.
\newblock Technical report, World Health Organization, Geneva, 2015.

\end{thebibliography}

\appendix
\section{Estimation of the rates from the survey}\label{app:estimation}

\emph{Data.} The survey took place in Tenkodogo, the main town of the Centre-Est region of Burkina Faso, in April 2025, in households of the formal (\emph{loti}) and informal (\emph{non loti}) zones that kept animals. In each household, a stool sample was taken from one member, faeces from one animal (poultry, small or large ruminants, donkeys or pigs) and soil from the yard. Samples were enriched in buffered peptone water and plated on MacConkey agar; \emph{Escherichia coli} was identified by cultural and biochemical tests, and its susceptibility to antibiotics was tested by disk diffusion \citep{bauer1966}. We use the results for five antibiotics, amoxicillin--clavulanic acid (AMC), tetracycline (TET), trimethoprim--sulfamethoxazole (SXT), gentamicin (GEN) and cefepime (FEP), in the $59$ households with usable results. Isolates were not available for every node of every household; the person and the animal were both tested in $20$ households.

\emph{Model and likelihood.} For each antibiotic $t$, the resistance of the person and of the animal of a household is described by the pair chain of Section~\ref{sec:model} with outside rates $\alpha_{H,t}$ and $\alpha_{A,t}$, a coupling $c$ common to the five antibiotics, and loss rates equal to $1$. The pair is assumed to be at its stationary law $\pi_t$, which is computed exactly from the $4\times4$ generator. A household in which both nodes were tested for antibiotic $t$ contributes $\log\pi_t(x_H,x_A)$ to the log-likelihood; a household in which only one node was tested contributes the logarithm of the marginal probability of its result. Since the five antibiotics are measured on the same isolates, the sum over antibiotics is a composite likelihood \citep{lindsay1988,varin2011}: it gives consistent estimates, but its curvature overstates the precision, so we do not use it for standard errors. We maximised it over the $11$ parameters with the L-BFGS-B method from eight random starting points.

\emph{Estimates.} The estimates for each antibiotic are given in Table~\ref{tab:estim}; the common coupling is $c=0.53$. We take as typical trait the geometric means of the outside rates over the five antibiotics, $\alpha_H=0.31$ and $\alpha_A=0.34$. The soil was tested in fewer households, and its prevalence, pooled over the five antibiotics, was $0.44$ in $70$ isolates; the rates involving the soil are therefore treated as assumptions (Section~\ref{sec:params}). The computation is in \texttt{estimate.py}.

\begin{table}[h]
\centering
\caption{Outside rates estimated for each antibiotic (time unit: mean duration of carriage).}
\label{tab:estim}
\small
\begin{tabular}{lcc}
\toprule
Antibiotic & $\alpha_{H,t}$ & $\alpha_{A,t}$\\
\midrule
AMC & 0.24 & 1.28\\
TET & 0.73 & 0.58\\
SXT & 0.71 & 0.18\\
GEN & 0.18 & 0.19\\
FEP & 0.12 & 0.18\\
\bottomrule
\end{tabular}
\end{table}

\section{Code}\label{app:code}
We checked Theorems~\ref{thm:comparison} and~\ref{thm:positive} numerically in \texttt{check\_theory.py}: for $60$ random parameter sets, random piecewise constant interventions $u\ge\tilde u$ and random initial laws $p_0\preceq\tilde p_0$ (with $\tilde p_0$ the law of the componentwise maximum of a state drawn from $p_0$ and an independent random state), the laws of household states were ordered on all $19$ non-empty increasing subsets of $\X$ at all times (smallest gap $-2\times10^{-15}$, that is zero up to rounding), and the adjoint was increasing in the household state at all times before the horizon (smallest increment $0.011$). The file \texttt{estimate.py} estimates the rates from the survey data (\texttt{data/households\_public.csv}, de-identified); \texttt{model.py} defines the household generator, the mean-field equation, the adjoint equation, the switching functions and the forward--backward sweep; \texttt{stoch.py} simulates $N$ households exactly; \texttt{run\_main.py}, \texttt{run\_sweeps.py}, \texttt{run\_2d.py} and \texttt{run\_mf.py} produce the results, \texttt{run\_synergy.py} and \texttt{synergy\_fixed.py} the study of synergy, \texttt{check\_numerics.py} the checks of the numerical method; \texttt{scalar.py} computes the function $g_u$ of Theorem~\ref{thm:scalar}, \texttt{model\_k.py} implements several types of households, \texttt{run\_revision.py} produces the results on equilibria, the long-run problem, the compartmental model, the two zones and the robustness study, and \texttt{figs.py} the figures.

\end{document}